\documentclass[a4paper,11pt]{article}
\usepackage{graphicx}
\usepackage[T1]{fontenc}
\usepackage[utf8]{inputenc}
\usepackage[english]{babel}
\usepackage{amsthm,amssymb,amsfonts}
\usepackage{amsmath}
\usepackage{lmodern}
\usepackage{subcaption}
\usepackage{geometry}
\usepackage[boxed]{algorithm2e}

\usepackage{framed}
\usepackage{color}
\definecolor{darkgreen}{rgb}{0.00,0.5,0.00}
\usepackage{hyperref}
\hypersetup{
	colorlinks=true,        
	linkcolor=blue,         
	citecolor=darkgreen,         
	urlcolor=blue           
}
\usepackage[backend=biber, style=numeric, 
doi=false, isbn=false, url=false,eprint=false,
giveninits=true,maxnames=6]{biblatex}
\usepackage{biblatex}

\usepackage{lscape} 
\usepackage{enumerate}

\newtheorem{theorem}{Theorem}[section]
\newtheorem{definition}{Definition}[section]
\newtheorem{lemma}{Lemma}[section]
\newtheorem{proposition}{Proposition}[section]
\newtheorem{remark}{Remark}[section]

\newtheorem{assumption}{Assumption}[section]

\newcommand{\R}{\mathbb R}
\newcommand{\N}{\mathbb N}
\renewcommand{\empty}{\varnothing}

\DeclareMathOperator{\prox}{prox}
\DeclareMathOperator{\sign}{sign}
\DeclareMathOperator{\proj}{proj}
\DeclareMathOperator{\acc}{acc}
\DeclareMathOperator{\conv}{conv}
\DeclareMathOperator{\clconv}{\overline{conv}}
\DeclareMathOperator*{\argmin}{argmin}
\DeclareMathOperator{\argstat}{stat-prox}
\DeclareMathOperator{\PhiD}{\Phi^{\rm DY}}
\DeclareMathOperator{\PhiP}{\Phi^{\rm PDHG}}
\DeclareMathOperator{\tPhiP}{\thickbar{\Phi}^{\rm PDHG}}
\DeclareMathOperator{\tPhiD}{\thickbar{\Phi}^{\rm DY}}

\usepackage{accents}
\newcommand\thickbar[1]{\accentset{\rule{.4em}{.8pt}}{#1}}

\newcommand{\diff}{\delta x}

\title{Splitting methods for nonconvex optimization in the absence of smoothness and weak convexity}
\author{Felipe Atenas\footnote{Centro de Modelamiento Matem\'{a}tico (CNRS IRL2807), Universidad de Chile, Santiago, Chile \newline(\href{mailto:fatenas@cmm.uchile.cl}{fatenas@cmm.uchile.cl})} \qquad Mikhail Solodov\footnote{IMPA -- Instituto de Matem\'atica Pura e Aplicada, Estrada Dona Castorina 110, Jardim Bot\^anico, Rio de Janeiro, RJ 22460-320, Brazil (\href{mailto:solodov@impa.br}{solodov@impa.br})}}
\date{\today}

\begin{document}
	\maketitle
	
	\begin{abstract}
		
	For nonconvex optimization problems, the existing convergence theory of proximal splitting 
methods depends on a trade-off between nonconvexity of some model components and smoothness or 
(strong) convexity in others. We show that such compensation is not necessary, 
at least for the
methods considered here. In fact, we do not  
 require even weak convexity for
any of the functions.
Assuming that each component of the objective function is only locally Lipschitz continuous and 
Clarke regular, we prove that bounded iterates generated by common proximal splitting methods, 
such as
generalizations of the Davis-Yin and the Primal-Dual Hybrid Gradient methods, 
accumulate at stationary points of the problem 
when the stepsizes vanish in the limit. 
Our proposal can also accommodate inertia/momentum terms in the iterates and 
computing merely
{\em stationary points} of proximal subproblems instead of proximal points prescribed by 
their global minimizers, which is a practically relevant distinction in nonconvex settings.
	
\end{abstract}
	
	\medskip
	
	
	\section{Introduction}

	We consider structured nonconvex optimization problems amenable to solve by means of 
proximal splitting methods. The first type of problems we study is the three-block optimization 
model	\begin{equation} \label{eq:DY-problem}
		\min_{x \in \R^d}\; f(x) + g(x) + h(x),
	\end{equation} where  $f, g, h: \R^d \to \R$ can all be nonconvex and
nonsmooth functions. Whenever $h$ is smooth, problem \eqref{eq:DY-problem} can be solved by the Davis-Yin splitting (DYS) method \cite{davis2017three},  which exploits the separable structure of the objective. The DYS method can be formulated as follows. Given an initial point $z^1 \in \R^d$ and a constant proximal/stepsize parameter $t >0$, 
the $k$-th iteration of DYS is defined as \begin{equation} \label{eq:DY}
		\left\{\begin{aligned}
			x^k & \in \prox_{t f}(z^k)\\
			y^k & \in \prox_{t g}(2x^k - z^k - t \nabla h(x^k))\\
			z^{k+1} & = z^k + y^k -  x^k,
		\end{aligned}\right.
	\end{equation} 
where $\prox$ denotes the usual proximal operator 
(the formal definition is in \eqref{eq:prox-argmin} further below). 
If $h=0$, the method in \eqref{eq:DY} reduces to the Douglas-Rachord splitting (DRS) method \cite{lions1979splitting}. In the convex setting, namely, when $f$, $g$, and $h$ are convex, if $h$ has Lipschitz continuous gradient with constant $\ell_h>0$, then the DYS sequences $\{x^k\}$ and $\{y^k\}$ asymptotically produce a global minimizer of \eqref{eq:DY-problem}, whenever $t \in (0, 2\ell_h^{-1})$. For nonconvex optimization problems, the state-of-the-art assumptions that guarantee subsequential convergence of the DYS sequences $\{x^k\}$ and $\{y^k\}$ to a stationary point
of \eqref{eq:DY-problem}, require both $f$ and $h$ to be smooth with Lipschitz continuous 
gradients, and $g$ may be nonsmooth nonconvex prox-bounded \cite{bian2021three}. 
Nevertheless, the DYS method does not make use of the gradient of $f$, 
and smoothness of $f$ thus appears as a trade-off for allowing nonconvexity of $g$. 
Setting $h=0$, the same smoothness assumption is needed to show convergence of the DRS method in nonconvex settings \cite{li2016douglas,themelis2020douglas}. 
In this work, we lift the smoothness assumptions on both $f$ and $h$, while still 
allowing $g$ to be nonconvex.

    The second type of problems we investigate are the composite models of the form \begin{equation} \label{eq:CP-problem-intro}
		\min_{x \in \R^d}\; f(Lx) + g(x) + h(x),
	\end{equation} where $f: \R^p \to \R$ and $g, h: \R^d \to \R$ can
 be nonconvex nonsmooth functions,  and $L: \R^d \to \R^p$ is a linear map. The problem in \eqref{eq:CP-problem-intro} can be solved, in principle, using DYS, in the following two settings: 
    \begin{enumerate}
        \item the functions $g$ and $h$ are ``prox-friendly'', and $f$ is smooth, or
        \item the functions $f \circ L$ and $g$ (resp. $h$) are ``prox-friendly'', and $h$ is smooth (resp. $g$).
    \end{enumerate}
The first possibility naturally fails whenever $f$ is nonsmooth. The second case depends on the nature of the operator $L$, as the computation of the proximal operator of the composition $f \circ L$ may not be a tractable operation.  As an alternative, we can solve the problem with the Primal Dual Hybrid Gradient (PDHG) method \cite{chambolle2011first,condat2013primal,vu2013splitting},
also known as the Chambolle-Pock method. 
If $h$ is smooth,  the method can be defined as follows. Given initial points 
$x^0 \in \R^d$, $y^0\in \R^p$, and $u^0 = x^0$,  constant proximal/stepsize 
parameters $\tau, \sigma > 0$, and a relaxation parameter $\theta \in [0,1]$, 
the $k$-th iteration of PDHG is \begin{equation} \label{eq:PDHG}
		\left\{\begin{aligned}
			y^{k+1} & \in \prox_{\sigma f^*}(y^k + \sigma Lu^k)\\
			x^{k+1} & \in \prox_{\tau g}(x^k - \tau L^*y^{k+1} - \tau \nabla h(x^k))\\
			u^{k+1} & = x^{k+1} + \theta (x^{k+1} - x^k).
		\end{aligned}\right.
	\end{equation} For convex problems, the original convergence results state that if $h$ has Lipschitz continuous gradient with constant $\ell_h>0$ and the stepsizes $\tau,\, \sigma >0$ satisfy $\frac{\tau \ell_h}{2} + \tau\sigma\|L\|^2 \leq 1$, then $\{(x^k,y^k)\}$ converges to a primal-dual solution pair of problem \eqref{eq:CP-problem-intro}, that is, a pair $(\bar x, \bar y)$ such that $0 \in L^* \bar y + \partial g(\bar x) + \partial h(\bar x)$ and $ L \bar x \in  \partial f^*(\bar y)$. In nonconvex settings, \cite{mollenhoff2015primal} shows convergence of $\{x^k\}$ to a unique global minimizer of problem \eqref{eq:CP-problem-intro} when $h=0$, $f$ is $\rho$-weakly convex, $g$ is $\mu$-strongly convex, such that $\mu > \rho \|L\|^2$. Under these assumptions, the objective function in \eqref{eq:CP-problem-intro} is thus strongly convex. In \cite{sun2018precompact}, subsequential convergence and global convergence is obtained assuming the Kurdyka--{\L}ojasiewicz property, allowing $g$ to be nonconvex nonsmooth, but $f$ needs to be smooth. In \cite{bolte2018nonconvex}, the authors introduce a Lagrangian methodology to solve problems generalizing \eqref{eq:CP-problem-intro} for $h=0$,  requiring smoothness of $g$ to compensate nonconvexity of $f$. 
 Further, \cite{bot2019proximal} analyzes a general framework, which when specialized to 
\eqref{eq:CP-problem-intro} with $h=0$ also needs a smooth-nonsmooth structure for convergence. 
Similarly as for the DYS method, we show that for the PDHG method,  
smoothness and strong convexity assumptions are not needed for allowing nonconvexity
 in the models.

    Since our work lifts structural assumptions in the existing literature, a new range of nonconvex optimization problems can be solved with theoretical guarantees of convergence. In particular, we only require the functions in the problems \eqref{eq:DY-problem} and \eqref{eq:CP-problem-intro} to be locally Lipschitz continuous and Clarke regular. In the following, we describe two structured classes of functions that fall into this category.

    First, consider the family of functions of the composite form $h = \rho \circ c$, where
$\rho : \R^m \to \R$ is a finite-valued convex Lipschitz continuous function, and
$c : \R^d \to \R^m$ is strictly differentiable. Such compositions are Clarke regular, and a
subgradient is obtained from the chain rule $\nabla c(x)^* s_\rho(c(x)) \in \partial h(x)$, where $s_\rho(c(x)) \in \partial \rho(c(x))$. The proximal operator of $h$ is rarely computable for such a composition. The analyses in \cite{bian2021three,alcantara2024four} require
$\nabla h$ to be Lipschitz continuous, and thus the outer function $\rho$ must be
replaced by a smooth surrogate, which perturbs the model. 

A second mechanism that yields nonconvex models outside the existing theory, even when $h$ is smooth, is having two nonsmooth proximable terms. The analyses in \cite{li2016douglas,themelis2020douglas,bian2021three} require $f$ to be differentiable with Lipschitz continuous gradient, so only one of the two proximable blocks is allowed to be nonsmooth. Our assumptions treat $f$ and $g$ symmetrically, and therefore cover any model in which two structurally different nonsmooth regularizers are imposed
simultaneously. This is typical for sparsity-promoting terms,  where a group
penalty and an entrywise penalty each admit a closed-form proximal operator while their sum does not.

Below, we describe three model instances that our framework covers, 
but which are outside the existing literature on splitting methods. 

\paragraph{Phase Retrieval.}
Let $a_1,\dots,a_m \in \R^d$ be measurement vectors and $b \in \R^m$ be
intensity measurements corrupted by outliers, and let $D \in \R^{p \times d}$
encode a discrete difference operator. Given a regularizer $p_2$  that induces sparsity in the solution, the model
\begin{equation} \label{eq:phase-retrieval}
	\min_{x \in \R^d} \; \underbrace{\lambda \|Dx\|_1}_{f(x)}
	\; + \; \underbrace{\sum_{i=1}^d p_2(x_i)}_{g(x)}
	\; + \; \underbrace{\frac{1}{m}\sum_{i=1}^m
		\left| \left<a_i,x\right>^2 - b_i \right|}_{h(x)}
\end{equation} fits our
 composite framework with $\rho = \frac{1}{m}\|\cdot\|_1$ and
$c(x) = (\left<a_i,x\right>^2 - b_i)_{i=1}^m$.
The misfit $h$ is the robust formulation of phase retrieval studied in \cite{duchi2018stochastic,davis2019stochastic}. This function $h$ is nonsmooth and
weakly convex on bounded sets, hence Clarke regular. Furthermore, $p_2$ can be taken to be, for instance, the $\ell_1$ norm, the \emph{smoothly clipped absolute deviation} \cite{fan2001variable}, and  the \emph{minimax concave penalty} \cite{zhang2010nearly}, all weakly convex, and hence Clarke regular, and the separable sum $g$ inherits this property. The existing analyses cannot 
accommodate this model due to the nonconvexity of $g$, and nonsmoothness of $h$ and $f$.  

\paragraph{Multiaffine Least Squares \cite{MR3418719}.}
Let $\mathcal{E}_1,\dots,\mathcal{E}_p$ and $\mathcal{F}$ be Euclidean spaces, and let $\mathcal{A}:\mathcal{E}_1\times\cdots\times\mathcal{E}_p\to\mathcal{F}$ be a multiaffine map, namely, affine in each block $x_i$ when the remaining blocks $j \neq i$ are fixed. Given $b\in\mathcal{F}$, a set $\mathcal{I}\subseteq\{1,\dots,p\}$ of blocks on which sparsity is desired, an entrywise regularizer $p_2$ as in \eqref{eq:phase-retrieval}, a smooth map $\mathcal{Q}$ with values in a Euclidean space, a norm $\nu$ on that space, and $\mu\geq0$, consider the model \begin{equation} \label{eq:multiaffine-model}
	\min_{x=(x_1,\dots,x_p)} \; \underbrace{\frac{1}{2}\|\mathcal{A}(x_1,\dots,x_p)-b\|^2}_{f(x)}
	\; + \; \underbrace{\sum_{i\in\mathcal{I}}\sum_{j} p_2\big([x_i]_j\big)}_{g(x)}
	\; + \; \underbrace{\mu\,\nu\big(\mathcal{Q}(x)\big)}_{h(x)},
\end{equation} where $[x_i]_j$ denotes the $j$-th entry of the block $x_i$. The block $f$ is the data-fitting term of a factorization model, $g$ promotes entrywise sparsity in the factors indexed by $\mathcal{I}$, and $h$ encodes structural requirements on the factors $x_1, \dots, x_p$. The latter is typically nondifferentiable at the points where these requirements hold, namely, $\mathcal{Q}(x)=0$. 

The block $f$ is polynomial, hence smooth, but may not be  weakly convex or have a globally Lipschitz continuous gradient. For example, one of the eigenvalues of the Hessian of $(u,v)\mapsto\frac{1}{2}(uv-m)^2$ at $u=v=a$ is $m-a^2$, which diverges as $|a|\to\infty$. Consequently, the analyses in \cite{bian2021three, li2016douglas, themelis2020douglas}, which require the smooth blocks to have a globally Lipschitz continuous gradient, do not apply to \eqref{eq:multiaffine-model}. 
However, all three blocks in \eqref{eq:multiaffine-model} are locally Lipschitz continuous and 
Clarke regular. In particular, the function $h$ is the composition of the finite convex function 
$\mu\nu$ with the smooth map $\mathcal{Q}$; hence, it is
 Clarke regular by the chain rule \cite[Theorem~2.6.6]{clarke1990optimization}.

Model \eqref{eq:multiaffine-model} covers the following instances.

\begin{enumerate}
	
	\item \emph{Robust PCA and matrix factorization.} Let $M\in\R^{m\times n}$ be observed on an index set $\Omega$, and suppose that the observed entries are the sum of a low-rank matrix and a sparse matrix of gross errors. Denote by $P_\Omega:\R^{m\times n}\to\R^{|\Omega|}$ the operator that extracts the entries indexed by $\Omega$. Take $x=(U,V,S)$ with factors $U\in\R^{m\times r}$ and $V\in\R^{n\times r}$ of the low-rank component, and corruptions $S\in\R^{|\Omega|}$ on the observed entries. Set $\mathcal{A}(U,V,S)=P_\Omega(UV^\top)+S$,  $b=P_\Omega(M)$ \cite{mateos2012robust}, and $\mathcal{I}=\{3\}$, so that $g$ promotes sparsity of the corruptions. Set $\mathcal{Q}(U,V,S)=U^\top U-V^\top V$ with $\nu=\|\cdot\|_F$ to remove the scaling ambiguity of the factors \cite{tu2016low} to distinguish the pair $(U,V)$ from $(UG, VG^{-\top})$, for $G$ nonsingular. 
Then \eqref{eq:multiaffine-model} models a robust PCA problem. 
Observe that taking $\mathcal{A}(U,V)=P_\Omega(UV^\top)$ recovers a \emph{matrix completion/factorization}  model \cite{burer2003nonlinear}.
	
	\item \emph{Dictionary learning.} To represent a matrix of observed data $Y\in\R^{n\times N}$, $K >n$, as a sparse linear combination of atoms (columns) from a dictionary $D\in\R^{n\times K}$ and a sparse code matrix $X\in\R^{K\times N}$, take $x=(D,X)$, and set $\mathcal{A}(D,X)=DX$, $b=Y$ and $\mathcal{I}=\{2\}$.  
Then  \eqref{eq:multiaffine-model} models the dictionary learning problem \cite{aharon2006k}. 
The map $\mathcal{Q}(D,X)=D^\top D-I_K$ with $\nu=\|\cdot\|_F$ penalizes atoms of non-unit norm and the overlap between any two distinct atoms (the so-called \emph{mutual coherence property}) of the dictionary \cite{yaghoobi2009parametric}.

\end{enumerate}

Our contributions can be summarized as follows. 
We lift the smoothness assumptions for convergence of proximal splitting methods in the 
nonconvex setting, 
allowing Clarke regular nonconvex (fully) nonsmooth optimization problems. 
In other words, our approach does not require a trade-off between nonconvexity of components
 with smoothness or strong convexity of other components.  
In particular, our analysis reinterprets DYS and PDHG as abstract subgradient methods 
without the need to define Lyapunov functions satisfying sufficient decrease conditions, 
which are typically defined ad hoc for each method in consideration.

Moreover, since convergence results in nonconvex optimization usually state that cluster points of 
the generated sequences are stationary points of the problem, it is only natural to weaken the 
requirement of computing global minimizers in the proximal steps in 
\eqref{eq:DY} and \eqref{eq:PDHG} to merely finding stationary points of the proximal 
minimization subproblems, an operation we define in \eqref{eq:argstat} below. 
 Furthermore, we also introduce momentum/inertia in the DYS method, 
a feature known to improve convergence of various algorithms.  
To the best of the authors' knowledge, this is the first work on proximal splitting methods for 
nonconvex nonsmooth optimization with this level of generality.

This work is organized as follows. In Section~\ref{s:prelim} we recall some basic 
variational analysis tools, and describe an abstract subgradient framework, which is
 the basis of our convergence analysis. In Section~\ref{s:stat-prox}, we define the \emph{stationary proximal point} operator, and discuss its relation to the classical
 proximal point operator. We then proceed to analyze in Section~\ref{s:n-DYS}
 a nonsmooth generalization of DYS 
 for solving general regular nonconvex structured optimization problems in the form 
of \eqref{eq:DY-problem}. We dedicate Section~\ref{s:PDHG} to the study of a variant of the PDHG method to solve problem \eqref{eq:CP-problem-intro} under nonconvex nonsmooth assumptions.

	\section{Notation and Preliminaries} \label{s:prelim}
	
	\subsection{Variational analysis tools}
	
We endow the space $\R^d$ with the standard Euclidean inner product and the induced norm. 
We denote by $\overline{X}$ the closure of a set $X \subseteq \R^d$, $\clconv X$ the closure
of the convex hull of $X$, and 
$\bar B(x,r)$ the closed ball centered at $x \in \R^d$ with radius $r > 0$. For 
a sequence $\{v^k\}$, we denote by $\acc(\{v^k\})$ the set of its accumulation points. 
We use $\proj_X$ to denote the Euclidean projection onto $X$.
In the product space $\R^d \times \R^d$, we denote by $\pi_1$ the projection operator onto 
the first space, namely, for $(x,u) \in \R^d \times \R^d$, we have $\pi_1(x,u) = x$. 

For a locally Lipschitz continuous function
 $F:\R^d\to \R$, the Clarke subdifferential of $F$ at 
$x \in\R^d$
  is given by 
\[ \partial F(x) = \conv\, \{ g\in \R ^d\mid
  \exists \, \{ y_k\} \subset {\cal D}_F:\,
   y_k \to x,\;  \nabla F(y_k) \to g\},\]
where ${\cal D}_F$ is the set of points at which $F$ is differentiable, and where $\conv D$ stands for the convex hull of a set $D$. The set $\partial F(x)$ is a nonempty closed convex set, and $\partial F$ is an upper semicontinuous, locally bounded mapping; see \cite{clarke1990optimization}. We denote by $S$ the set of stationary points of $F$, that is, 
$$S = \{x \in \R^d: 0 \in \partial F(x)\}.$$ 
The function $F$ is Clarke regular at $x$ if $\partial F (x)$ coincides with the 
limiting subdifferential \cite[Corollary 8.11]{rockafellar2009variational}. 
We say that $F$ is $r$-weakly convex for some $r \geq 0$ if $F(\cdot) + \frac{r}{2}\|\cdot\|^2$ is convex. In particular, weakly convex functions are Clarke regular at any point. Given $y \in \R^d$, the proximal operator of a function $F$ is given by \begin{equation} \label{eq:prox-argmin}
		\prox_{\phi}(y) = \argmin_{w\in\R^d} \left\{F(w) + \frac{1}{2}\|w-y\|^2\right\}.\end{equation} We say $F$ is prox-bounded if there exists $t>0$ such that $F(w) + \frac{1}{2t}\|w-y\|^2$ is bounded below. In this case, $\prox_{tF}(y)$ is nonempty.  If $F$ is a $r$-weakly convex function, then $\prox_{tF}(y)$ is well-defined and single-valued for all $t \in (0, r^{-1})$ \cite[Proposition 12.19]{rockafellar2009variational}.

\subsection{Abstract subgradient framework}\label{sec:abstract}

We now describe a framework for analyzing convergence
of abstract subgradient-related iterative processes \cite{zavp:90a};
see also \cite{solzav:gr}.

Let $V \subseteq \R^d$ be an open set, and let $\{v^k\}$ be a sequence in $\R^d$ such that \begin{equation} \label{eq:gen-bounded-assumption}
	\{v^k\}\text{~is bounded, and for some compact convex set $C \subseteq V$,~} \acc(\{v^k\}) \subseteq C.
\end{equation} Let $\mathcal{L}: \R^d \to \R$ be a Clarke-regular, locally Lipschitz continuous 
function. For $v \in \R^d$, define the set 
\begin{equation} \label{eq:R}
	R(v) = \conv\{\partial \mathcal{L}(v) \cup N_C(v)\},
\end{equation}
where $N_C(v)$ is the standard normal cone to the convex set $C$, i.e.,
it is empty if $v\not\in C$, and for $v\in C$ it is given by
$\{z\in \R^d : \left<z,y-v\right> \leq 0 \;\forall\, y\in C\}$.
 
 For  a set-valued map $\Phi: \N \times \R^d \to 2^V$, which defines the iterative
process (see \eqref{iter} below), define 
 \begin{equation*}
\thickbar{\Phi}(v) := \limsup_{\substack{k \to \infty \\ \hat{v} \to v}} \Phi(k,\hat{v}),
\end{equation*} 
the outer limit of $\Phi$ at $v$. Furthermore, define the set \begin{equation} \label{eq:A}
	A(\mathcal{L},C,\Phi) = \left\{v \in C: \max_{r \in R(v)}\min_{\phi\in\thickbar{\Phi}(v)} \left<r,\phi\right> \leq 0\right\}.
\end{equation} 
According to Theorem~\ref{t:subgrad-scheme-convergence} stated below,
the set $A(\mathcal{L},C,\Phi)$ serves as the attractor for
the iterates generated according to the map $\Phi$. 

A set $\Omega \subseteq \R^d$ is said to be $\mathcal{L}$-connected if the image $\mathcal{L}(\Omega) = \{ t \in \R: \exists v \in \Omega, t = \mathcal{L}(v)  \}$ is a connected set in $\R$. Consider a decomposition $\{A^j\}_{j \in J}$ of $	A(\mathcal{L},C,\Phi)$ into its $\mathcal{L}$-connected components, namely, \begin{equation}\label{eq:connected-components}
		A(\mathcal{L},C,\Phi) = \bigcup_{j \in J} A^j, \quad
 A^{j_1} \neq A^{j_2}, \, j_1 \neq j_2, \, j_1, j_2 \in J.
\end{equation}

The following result describes  
convergence properties of an abstract subgradient scheme associated with the set-valued map $\Phi$.
It is as given in {\cite[Theorem 3.1]{solodov2025convergence}}, 
which is a simplified version of the statements in 
\cite{zavp:90a},  \cite{solzav:gr}.
 
\begin{theorem} \label{t:subgrad-scheme-convergence}
	Let $V \subseteq \R^d$ be an open set, and let $\Phi : \N\times V \to 2^V$ be a set-valued map. For $v^1 \in V$ and a positive sequence of parameters $\{\alpha_k\}$, consider the sequence $\{v^k\}$  in $\R^d$ generated via 
\begin{equation}\label{iter}
		v^{k+1} \in v^k - \alpha_k \Phi(k, v^k) \text{~for all~} k \in \N,
\end{equation} 
such that \eqref{eq:gen-bounded-assumption} holds. Suppose further that the following conditions hold: \begin{equation} \label{t:subgrad-scheme-convergence-1}
	\sup_{\hat{v}\in V} \limsup_{\substack{v \to \hat{v}\\ k \to \infty}} \sup_{\phi \in \Phi(k,\hat{v})}\|\phi\|<\infty,
	\end{equation} and 	
	\begin{equation} \label{t:subgrad-scheme-convergence-2}
		\alpha_k \to 0, \qquad \sum_{k=1}^\infty \alpha_k = \infty.
	\end{equation}  Then,  there exists $j \in J$ from \eqref{eq:connected-components} such that all the accumulation points of the sequence $\{\mathcal{L}(v^k)\}$ belong to $\mathcal{L}( A^j \cap \acc(\{v^k\}))$. In addition, if the set $\mathcal{L}(A(\mathcal{L},C,\Phi))$ is nowhere dense in $\R$, then $\acc(\{v^k\})$ is a subset of an $\mathcal{L}$-connected component of $A(\mathcal{L},C,\Phi)$.
\end{theorem}

\section{The stationary proximal map}\label{s:stat-prox}

 The proximal point operator is the cornerstone of a plethora of methods in nonsmooth optimization. However, its definition (see \eqref{eq:prox-argmin}) requires a global minimizer of the proximal minimization subproblem. For general nonconvex functions, computing a proximal point may be restrictive, as general optimization solvers can only guarantee, in principle, convergence to stationary points. This motivates the definition of the \emph{stationary proximal point operator}.

\begin{definition}[Stationary proximal point operator] \label{d:stat-prox}
	Given a function $\phi: \R^d \to \R $, the stationary proximal map $\argstat_\phi: \R^d \rightrightarrows \R^d$ is defined as \begin{equation} \label{eq:argstat}
		x \in \argstat_\phi(y) \iff 0 \in \partial \phi(x) + x - y, \quad \text{for all~} x,y \in \R^d.
	\end{equation} 
\end{definition}

In other words, $\argstat_\phi$ is the map of stationary points associated with the proximal minimization problem of $\phi$, namely, \begin{equation} \label{eq:prox-minimization}
	\min_{w\in\R^d} \left\{\phi(w) + \frac{1}{2}\|w-y\|^2\right\}.
	\end{equation} 
Clearly, for all $y\in\R^d$ and $t>0$, \begin{equation} \label{eq:prox-subset-argstat}
	\prox_{t\phi}(y)\subseteq\argstat_{t\phi}(y).
\end{equation} In particular, if $\phi$ is $r$-weakly convex for some $r \geq 0$, then $\prox_{t\phi}(y) = \argstat_{t\phi}(y)$ is a singleton for all $t \in (0,r^{-1})$.

    In the next result, we first characterize the existence of stationary proximal points for locally Lipschitz functions.

    \begin{lemma}[The stationary proximal map is well-defined]\label{l:argstat-existence}
Let $\phi:\R^d\to\R$ be locally Lipschitz continuous, $y\in\R^d$ and $r>0$, and set
\[
L_r(y) := \sup\left\{\|s\| : s\in\partial\phi(w),\ w \in \bar B(y,r)\right\} < \infty .
\]
If $t\in (0, rL_r(y)^{-1}]$, then there exists $x\in\argstat_{t\phi}(y) \cap \bar B(y,t\,L_r(y))$.
In particular, if $\phi$ is globally $L$-Lipschitz continuous, then for all $t>0$ and $y\in\R^d$, $\argstat_{t\phi}(y)\neq\empty$, and every $x\in\argstat_{t\phi}(y)$ satisfies $x \in \bar B(y, tL)$.
\end{lemma}

\begin{proof}
Let $K:=\bar B(y,r)$ and $T(w):=y-t\,\partial\phi(w)$ for $w\in K$. Since $\partial\phi$ has nonempty convex compact values, closed graph, and is bounded on $K$, the map $T$ is upper semicontinuous with nonempty convex compact values, and $T(K)\subseteq \bar B(y,tL_r(y))\subseteq K$. By Kakutani's fixed point theorem \cite{Kakutani1941AGO}, there exists $x\in K$ with $x\in T(x)$, that is, $0\in t\partial\phi(x)+x-y$, such that $\|x-y\|\leq tL_r(y)$. In the globally Lipschitz case, take $r=tL$. The bound $\|x-y\|=t\|s\|\leq tL$ holds for every stationary point since $s\in\partial\phi(x)$.\end{proof}

In view of \eqref{eq:prox-subset-argstat}, the stationary proximal map enters the picture when there exist stationary points of \eqref{eq:prox-minimization} that are not global minimizers, or when the proximal operator produces empty values. In the next result, we provide an explicit threshold for 
the proximal parameter that determines whether the proximal operator is empty-valued.
 
\begin{lemma} \label{l:prox-dichotomy}
	Let $\phi:\R^d\to\R$ be a lower semicontinuous prox-bounded function and $t>0$. Let \begin{equation} \label{eq:decay-constant}
	c_\phi := \max\left\{0,\; -\liminf_{\|x\|\to\infty}\frac{2\phi(x)}{\|x\|^2}\right\}\in[0,+\infty].
\end{equation} Then, the following holds.
	\begin{enumerate}[(i)]
		\item \label{l:prox-dichotomy-i} If $t\,c_\phi<1$, then $\prox_{t\phi}(y)\neq\emptyset$ for all $y\in\R^d$.
		\item \label{l:prox-dichotomy-ii} If $t\,c_\phi>1$, then $\prox_{t\phi}(y)=\emptyset$ for all $y\in\R^d$.
		\item \label{l:prox-dichotomy-iii} If $\phi$ is $r$-weakly convex, then $c_\phi\leq r$. If $\phi$ is bounded from below, then $c_\phi=0$.
	\end{enumerate}
\end{lemma}
 
\begin{proof}
	From \cite[Exercise 1.24]{rockafellar2009variational}, $\frac{1}{c_\phi}$ is the threshold of prox-boundedness of $\phi$, and thus item \eqref{l:prox-dichotomy-i} follows. For item \eqref{l:prox-dichotomy-ii}, if $t\,c_\phi>1$, the problem in \eqref{eq:prox-minimization} is unbounded, and thus $\prox_{t\phi}(y)=\emptyset$. For item \eqref{l:prox-dichotomy-iii}, for any $v \in \partial \phi(0)$, we have $\phi(x)\geq \phi(0)+\left<v,x\right>-\frac{r}{2}\|x\|^2$ for all $x\in\R^d$, which gives $\liminf_{\|x\|\to\infty}2\phi(x)/\|x\|^2\geq-r$, that is, $c_\phi\leq r$. If $\phi\geq m$ for some $m\in\R$, then $\liminf_{\|x\|\to\infty}2\phi(x)/\|x\|^2\geq0$, so $c_\phi=0$.\end{proof}
  
Observe that when $t c_\phi= 1$, it is not possible, in general, 
to determine whether $\prox_{t \phi}$ is nonempty. 
	
For weakly convex functions, the proximal point operator is well-defined for small enough 
proximal parameters, but it may be empty-valued for large parameters. Table~\ref{tab:regimes} summarizes the relationship between the proximal operator and the stationary proximal operator of an $r$-weakly convex function as a function of the proximal parameter $t$, the constant of weak convexity $r$ and the asymptotic decay constant $c_\phi$.
 
 \begin{table}[!h]
 	\centering
 	\renewcommand{\arraystretch}{1.4}
 	\begin{tabular}{|c|p{4.3cm}|p{5.6cm}|}
 		\hline
 		 Parameter bounds & $\prox_{t\phi}(y)$ & $\argstat_{t\phi}(y)$ \\ \hline
 		 $0<t<1/r$ & singleton & coincides with $\prox_{t\phi}(y)$ \\ \hline
 		$1/r\leq t<1/c_\phi$ & nonempty and possibly multivalued & contains $\prox_{t\phi}(y)$, possibly strictly\\ 
 		\hline
 		 $t>1/c_\phi$ & empty & possibly nonempty (Lemma~\ref{l:argstat-existence}) and no point is a global minimizer of $t\phi(\cdot) + \frac{1}{2}\|\cdot-y\|^2$ \\ \hline
 	\end{tabular}
 	\caption{Regimes of $\prox_{t\phi}$ and $\argstat_{t\phi}$ for a $r$-weakly convex function $\phi$ with constant $c_\phi\leq r$ defined in \eqref{eq:decay-constant}.}
 	\label{tab:regimes}
 \end{table}
 
 The first regime in Table~\ref{tab:regimes} is the classical parameter setting for weakly convex functions. However, it may impose restrictive bound constraints for the proximal parameter. The convergence of proximal-based methods may benefit from taking a large proximal parameter, although the proximal operator may not be well-defined in this case. Hence, the stationary proximal map appears as an alternative to bridge this gap.

 Table~\ref{tab:regimes} suggests that even for weakly convex functions, the stationary proximal point map may be strictly larger that the proximal point operator. We next illustrate 
this phenomenon for the the minimax concave penalty \cite{zhang2010nearly}. 
 This function is used to induce sparsity in variable selection problems, and reduces the \emph{shrinking bias} of the usual $\ell_1$ penalty, see \cite{MR4878992} and references therein. 
 We show that for the minimax concave penalty we can obtain closed-form formulae for the stationary proximal point operator similar to those of its  (pure) proximal point map, including the limiting case $t c_\phi = 1$.

 \begin{lemma}[$\argstat$ of minimax concave penalty] \label{l:stat-prox-MCP}
 	Given $\lambda,\,\gamma>0$ and $w \in \R$. Consider the minimax concave penalty given by 
 	\begin{equation*}
 		\phi(w)=\begin{cases}
 			\lambda|w|-\dfrac{w^2}{2\gamma}, & |w|\leq\gamma\lambda,\\
 			\dfrac{\gamma\lambda^2}{2}, & |w|>\gamma\lambda.
 		\end{cases}
 	\end{equation*}
Then $\argstat_{t \phi}(y) \neq \emptyset$ for all $y \in \R$ and $t >0$. Furthermore,	if $t<\gamma$, 
 	\begin{equation*}
 		\argstat_{t\phi}(y)=\prox_{t\phi}(y)=\begin{cases}
 			\{0\}, & |y| \leq t\lambda,\\
 			\left\{\sign(y)\dfrac{\gamma(|y|-t\lambda)}{\gamma-t}\right\}, & t\lambda<|y|<\gamma\lambda,\\
 			\{y\}, & |y|\geq\gamma\lambda.
 		\end{cases}
 	\end{equation*}
 	If $t=\gamma$, then 
 	\begin{equation} \label{eq:argstat-MCP-case-2}
 		\prox_{t\phi}(y)=\argstat_{t\phi}(y)=\begin{cases}
 			\{0\}, & |y| < \gamma \lambda,\\
 			\sign(y)[0,\gamma\lambda], & |y|=\gamma\lambda,\\
 			\{y\}, & |y|>\gamma\lambda.
 		\end{cases}
 	\end{equation}
 	 If $t>\gamma$, then
 	\begin{equation} \label{eq:argstat-MCP-case-3}
 		\argstat_{t\phi}(y)=\begin{cases}
 			\{0\}, & |y|<\gamma\lambda,\\
 			\{0,\sign(y)\gamma\lambda\}, & |y|=\gamma\lambda,\\
 			\left\{0,\sign(y)\dfrac{\gamma(t\lambda-|y|)}{t-\gamma},y\right\}, & \gamma\lambda<|y|<t\lambda,\\
 			\{0,\sign(y)t\lambda\}, & |y|=t\lambda,\\
 			\{y\}, & |y|>t\lambda,
 		\end{cases} 
 	\end{equation} and \begin{equation} \label{eq:prox-MCP-case-3-prox}
 	\prox_{t\phi}(y)=\begin{cases}
 		\{0\}, & |y|<\lambda\sqrt{t\gamma},\\
 		\{0,y\}, & |y|=\lambda\sqrt{t\gamma},\\
 		\{y\}, & |y|>\lambda\sqrt{t\gamma}.
 	\end{cases}
 	\end{equation} In particular, $\argstat_{t\phi}(y)\supsetneq\prox_{t\phi}(y)$ if and only if $\gamma\lambda\le|y|\le t\lambda$.
 \end{lemma}
 
 \begin{proof}
 	The function $\phi$ is continuously differentiable on $\R\setminus\{0\}$ with
 	\begin{equation} \label{eq:MCP-derivative}
 		\phi'(w)=\sign(w)\left(\lambda-\frac{|w|}{\gamma}\right) \text{ for } 0<|w|\leq\gamma\lambda, \quad \phi'(w)=0 \text{ for } |w|\geq\gamma\lambda,
 	\end{equation}
 	and, at the origin, $\partial\phi(0)=[-\lambda,\lambda]$. Hence, $\partial \phi(w) \subseteq [-\lambda,\lambda] $ for all $w \in \R$, and thus $\phi$ is $\lambda$-Lipschitz continuous \cite[Theorem 9.13]{rockafellar2009variational}. Lemma~\ref{l:argstat-existence} yields $\argstat_{t \phi}(y) \neq \emptyset$ for any $y \in \R$ and $t >0$.
 	
 	  Moreover,  $\phi$ is $r$-weakly convex with $r=1/\gamma$ \cite{zhang2010nearly}. Since $\phi\geq0$, $c_\phi=0$ by Lemma~\ref{l:prox-dichotomy}\eqref{l:prox-dichotomy-iii}. Therefore, the first regime in Table~\ref{tab:regimes} corresponds to $t<\gamma$ (the usual weakly convex proximal parameter bound), and thus $\argstat_{t \phi}$ is the \emph{firm thresholding} operator \cite{MR1488646}. The second case in Table~\ref{tab:regimes} is  $t\geq\gamma$, and the third case does not apply. Since $\phi$ is even, it follows that $\partial\phi(-w)=-\partial\phi(w)$ for all $w\in\R$, so that $\argstat_{t\phi}(-y)=-\argstat_{t\phi}(y)$ for all $y\in\R$.  It thus suffices to compute $\argstat_{t\phi}(y)$ for $y\geq0$. By \eqref{eq:argstat}, $w\in\argstat_{t\phi}(y)$ if and only if $y-w\in t\,\partial\phi(w)$. We distinguish four cases depending on the value of $w$.
 	\begin{enumerate}[(a)]
 		\item If $w=0$, then the stationarity condition reads $y\in[-t\lambda,t\lambda]$, so that $y\leq t\lambda$.
 		
 		\item If $0<w<\gamma\lambda$, then, from \eqref{eq:MCP-derivative}, the condition becomes $y-w=t\lambda-\frac{t}{\gamma}w$, that is,
 		\begin{equation} \label{eq:MCP-middle}
 			\left(1-\frac{t}{\gamma}\right)w=y-t\lambda.
 		\end{equation}
 		(b.1) If $t<\gamma$, then \eqref{eq:MCP-middle} gives $w=\frac{\gamma(y-t\lambda)}{\gamma-t}$, and the constraint $0<w<\gamma\lambda$ is equivalent to $t\lambda<y<\gamma\lambda$.
 		(b.2) If $t>\gamma$, then \eqref{eq:MCP-middle} gives $w=\frac{\gamma(t\lambda-y)}{t-\gamma}$. Since $t-\gamma>0$, the inequality $\frac{\gamma(t\lambda-y)}{t-\gamma}>0$ holds if and only if $y<t\lambda$. Moreover, $\frac{\gamma(t\lambda-y)}{t-\gamma}<\gamma\lambda$ if and only if $t\lambda-y<\lambda(t-\gamma)$, that is, $y>\gamma\lambda$. Hence, the constraint $0<w<\gamma\lambda$ is equivalent to $\gamma\lambda<y<t\lambda$.
 		(b.3) If $t=\gamma$, then \eqref{eq:MCP-middle} gives $0=y-\gamma\lambda$. Thus, every $w\in(0,\gamma\lambda)$ is stationary when $y=\gamma\lambda$, and no $w$ is stationary otherwise.
 		
 		\item If $w\geq\gamma\lambda$, then $\partial\phi(w)=\{0\}$ by \eqref{eq:MCP-derivative}, so the stationarity condition is $w=y$, which is admissible if and only if $y\geq\gamma\lambda$.
 		
 		\item If $w<0$, for $-\gamma\lambda<w<0$, the stationarity condition is $y=\left(1-\frac{t}{\gamma}\right)w-t\lambda$. If $t\leq\gamma$, then $y \leq -t\lambda<0$, which is not possible due to the assumption $y\geq 0$. If $t>\gamma$, then $\left(1-\frac{t}{\gamma}\right)w-t\lambda = \left(\frac{t}{\gamma}-1\right)|w|-t\lambda< \left(\frac{t}{\gamma}-1\right)\gamma\lambda-t\lambda=(t-\gamma)\lambda-t\lambda=-\gamma\lambda<0$, which is again incompatible with the assumption. Next, suppose $w\leq-\gamma\lambda$, so the stationarity condition implies $y=w<0$, a contradiction as well. In this manner, in all cases $y<0$, so there are no negative stationary points when $y\geq0$.
 	\end{enumerate}
 	
 	Hence, we have the following for $\argstat_{t \phi}(y)$ and $y\geq 0$. If $t=\gamma$, we recover the first line in \eqref{eq:argstat-MCP-case-2} from case (a), the second line from cases (a) and (b.3), and the third line from case (c). Moreover, since $\phi$ is $(1/\gamma)$-weakly convex, it follows that $\phi(\cdot) + \frac{1}{2t}\|\cdot - y\|^2$ is convex (not strongly convex), so that any stationary point is a global minimizer. That is, in this case, any stationary proximal point is a proximal point, which completes the proof of the identities in \eqref{eq:argstat-MCP-case-2}.
 	
 	If $t>\gamma$, then we recover the first line in  \eqref{eq:argstat-MCP-case-3} from case (a), the second line from cases (a) and (c), and third line from cases (c) and (b.2), the fourth line from cases (a) and (c), and the fifth line from case (c). We now compute $\prox_{t \phi}$. Given $y \geq 0$, let $\phi_y(\cdot) = t \phi(\cdot) + \frac{1}{2}(\cdot - y)^2$, so that \begin{equation*}
 	\phi_y(w)=\begin{cases}
 		t\lambda|w|-\dfrac{1}{2}\left(\dfrac{t}{\gamma}-1\right)w^2 - wy + \frac{y^2}{2}, & |w|\leq\gamma\lambda,\\
 		\dfrac{t\gamma\lambda^2}{2} + \dfrac{1}{2}(w-y)^2, & |w|>\gamma\lambda.
 	\end{cases}
 	\end{equation*} Then $\phi_y$ is coercive, and thus has at least one global minimizer contained in $\argstat_{t\phi}$. Hence, choosing the stationary point(s) in \eqref{eq:argstat-MCP-case-3} with smallest value evaluated at $\phi_y$ yields \eqref{eq:prox-MCP-case-3-prox}.\end{proof}

\begin{remark}[Approximate stationary proximal points] \label{r:argstat-multiaffine}
	Definition~\ref{d:stat-prox} is the exact version of the stationary proximal point operator. In practice, and similarly to the proximal point operator, the calculation of a stationary proximal point may require a subroutine which can only guarantee the computation of an \emph{approximate} stationary proximal point, defined as follows. Given $\varepsilon>0$ and $y\in \R^d$, let \begin{equation*}
		x = \argstat^\varepsilon_{\phi}(y) \quad\text{such that}\quad 0\in\partial\phi(x)+x-y+\varepsilon\bar B(0,1),
	\end{equation*}
	in the spirit of \cite{solodov2025convergence}. This notion is useful in the multiaffine least squares problem of \eqref{eq:multiaffine-model} described in the Introduction. More precisely, let $f$ be the data-fitting block in  \eqref{eq:multiaffine-model}, let $t>0$ and let $y=(y_1,\dots,y_p)$. Computing $\prox_{tf}(y)$ demands a global minimizer of the nonconvex polynomial $\psi_y(\cdot) := tf(\cdot)+\frac{1}{2}\|\cdot-y\|^2$, which is hard to find in general. In contrast, applying an optimization subroutine to $\psi_y$, such as the proximal alternating linearized minimization method (PALM) \cite{MR3232623}, generates stationary points of $\psi_y$ in the limit, that is, elements of $\argstat_{tf}(y)$. In practice, the subroutine exits before a 
stationary point is reached or, equivalently, upon the calculation of an approximate stationary proximal point. The statements of Algorithms~\ref{a:DY-shadow} and~\ref{a:PDHG} below use the stationary proximal point of Definition~\ref{d:stat-prox}, but it can be substituted by its approximate counterpart.
 \end{remark}

\section{Nonconvex DYS method without smoothness} \label{s:n-DYS}

    In this section, we introduce a three-block proximal splitting method that generalizes 
the DYS method to find stationary points of the objective function of \eqref{eq:DY-problem}. 
Given locally Lipschitz continuous functions $f,\,g, \, h: \R^d \to \R$, 
our goal is to find a point $\bar x \in \R^d$ such that \begin{equation*}
    	0 \in \partial f(\bar x) + \partial g(\bar x) + \partial h(\bar x).
    \end{equation*}

To motivate the method we study in this section, we first reformulate DYS to make the 
so-called shadow sequence $\{x^k\}$ the ``leading'' update, an idea inspired 
by \cite{atenas2025shadow}. 
Let $x^1 \in \R^d$, $t >0$, and set $z^1 = x^1 + ts_f(x^1)$ for some $s_f(x^1) \in \partial f(x^1)$. Then, $x^1 \in \argstat_{tf}(z^1)$. Furthermore, for all $k\in\N$, define the iterates via \begin{equation} \label{eq:DYb}
		\left\{\begin{aligned}
			y^k & \in \argstat_{t g}(2x^k - z^k - \gamma s_h(x^k)), \quad \mbox{for some } s_h(x^k) \in \partial h(x^k)\\
			z^{k+1} & = z^k + y^k -  x^k\\
			x^{k+1} & \in \argstat_{t f}(z^{k+1}).
		\end{aligned}\right.
	\end{equation} If $h$ is smooth, changing the order of the update rules, and substituting the proximal steps with stationary proximal steps, the schemes \eqref{eq:DYb} and \eqref{eq:DY} generate the same sequences $\{x^k\}$, $\{y^k\}$ and $\{z^k\}$, assuming
 the same selection of (stationary) proximal points. 
The reformulation in \eqref{eq:DYb} enables a novel interpretation of DYS as an abstract subgradient scheme for the so-called \emph{shadow}  sequence $\{x^k\}$. In contrast, common approaches require the use of the fixed-point sequence $\{z^k\}$ to lead the analysis.

	In addition to the reordering above, we also add a momentum term to $\{x^k\}$, 
the ``leading''  sequence that will actually produce stationary points asymptotically. 
Algorithm~\ref{a:DY-shadow} is an inertial version of the DYS method without smoothness assumptions, in which the gradient step on $h$ is replaced by a subgradient step, and both proximal steps are replaced by stationary proximal steps. 
Our convergence results apply to \emph{any} selection from $\argstat$, 
even when no proximal point exists (the corresponding subproblem has no global minimizers).

	\begin{algorithm}[!h]
		\caption{Three-block splitting with momentum for solving problem \eqref{eq:DY-problem}. \label{a:DY-shadow}}
		\SetKwInOut{Input}{Input}
		\Input{Choose starting points $x^0, x^1 \in \R^d$. 
			Choose $s_{f}(x^1) \in \partial f(x^1)$ and set $z^{1} = x^1 + \tau_k^1 s_{f_1}(x^1)$.}
		
		\For{$k=1,2,\dots$}{
			Step 1. Choose $\sigma_k, \, \alpha_k>0$, and compute $s_{h}(x^{k}) \in \partial h(x^k)$ and $$y^{k} \in \argstat_{\sigma_k g}(2x^{k} - z^{k} - \alpha_k s_{h}(x^{k})).$$
			
			Step 2. Update $z^{k+1} = z^{k} + y^{k} - x^{k}$.
			
			Step 3. Choose $\tau_k>0$ and compute $u^{k} \in \argstat_{\tau_k f}(z^{k+1})$.
			
			Step 4. Choose $\gamma_k >0$ and update $x^{k+1} =u^{k} + \gamma_k(x^k - x^{k-1})$.
		}
	\end{algorithm}

    We work under the following blanket assumptions to investigate 
 convergence properties of Algorithm~\ref{a:DY-shadow}.
    
     \begin{assumption} \label{as:DY}
    	For Algorithm~\ref{a:DY-shadow}, suppose the following conditions hold.
    	\begin{enumerate}[(i)]
    		\item \label{as:DY-i} The functions $f,\, g, \, h: \R^d \to \R$  are locally Lipschitz continuous and Clarke regular. 
    		\item \label{as:DYii} The parameters of Algorithm~\ref{a:DY-shadow} satisfy
    		\begin{equation} \label{eq:stepsize-assumption}
    			\frac{\tau_k}{\alpha_k} \to 1, \quad \frac{\sigma_k}{\alpha_k} \to 1,\quad 
\frac{\gamma_k}{\alpha_k} \to 0, \quad \frac{\gamma_{k-1}}{\alpha_k} \to 0, \quad
\alpha_k \to 0\; \text{as~} k \to \infty, \; \text{and~} 
 \sum_{k=1}^\infty \alpha_k = \infty .
    		\end{equation}
    		
    	\end{enumerate}
    \end{assumption}

	We shall show that under Assumption~\ref{as:DY}, 
the cluster points of the primal sequences $\{x^k\}$ and $\{u^k\}$ generated by 
Algorithm~\ref{a:DY-shadow} are stationary points for \eqref{eq:DY-problem}.
	
		\subsection{The DYS abstract subgradient scheme} \label{ss:scheme-DY}

To establish convergence properties of Algorithm~\ref{a:DY-shadow}, 
we first cast its iterates as abstract subgradient steps 
in the framework of Section~\ref{sec:abstract}.
	
	\begin{proposition} \label{p:DY-shadow-as-subgradient-step}
		
		Let Assumption~\ref{as:DY} hold. Suppose the sequences $\{x^k\}$, $\{u^k\}$ and $\{y^k\}$ generated by Algorithm~\ref{a:DY-shadow} are bounded, and consider an open bounded 
convex set $X \subseteq \R^d$ containing these sequences. For a constant $c>0$ dependent on $X$, let the set-valued mapping $\PhiD: \N\times X \times X \to 2^{X} \times 2^{X} $ be given by, for all $x, u\in \R^d$,  \begin{equation} \label{eq:def-Phi}
		(\phi_x, \phi_u) \in \PhiD(k,x,u) \iff\left\{\begin{aligned}
			&\phi_x = d - \frac{\gamma_k}{\alpha_k}\diff^+, \quad \phi_u = d - \frac{\gamma_{k-1}}{\alpha_k}\diff, \\ &d = \frac{\tau_k}{\alpha_k}s_f(q)+ \frac{\sigma_k}{\alpha_k}s_{g}(y) +  s_h(x), \\
			& s_f(q) \in \partial f(q),\, s_{g}(y) \in \partial g(y), \,  s_h(x) \in \partial h(x),\\
			& q = u - \alpha_k \left(d - \frac{\gamma_{k-1}}{\alpha_k}\diff\right), \\
			& y = 2x - u - \alpha_k s_h(x) - \sigma_k s_g(y) - \tau_{k-1} s_f(u), \, \|y\| \leq c, \\
			& \diff^+, \diff \in \R^d, \, \|\diff^+\| \leq c, \|\diff\| \leq c.
		\end{aligned}\right.
	\end{equation}
	 Then, \begin{equation*} 
	\begin{pmatrix}
		x^{k+1} \\ u^k
	\end{pmatrix} \in \begin{pmatrix}
		x^{k} \\ u^{k-1}
	\end{pmatrix} - \alpha_k \PhiD(k,x^k,u^{k-1}) \quad \text{for all } k \in \N.
	\end{equation*} 
		
	\end{proposition}

	\begin{proof}
	In view of Steps~$1$ and~$2$ in Algorithm~\ref{a:DY-shadow}, there exist $s_{g}(y^{k}) \in \partial g(y^{k})$ and $s_{f}(u^{k}) \in \partial f(u^{k})$ such that \begin{align}
		0 & = \sigma_k s_{g}(y^{k}) + y^{k} - (2x^{k} - z^{k} - \alpha_k s_{h}(x^{k})), \label{eq:OC-v} \\
		0 & = \tau_k s_{f}(u^{k}) + u^{k} - z^{k+1}. \label{eq:OC-u}
	\end{align} 
Combine \eqref{eq:OC-v}, and \eqref{eq:OC-u} with $k \gets k-1$, to obtain \begin{align*}
	y^k & = 2x^{k} - z^{k} - \alpha_k s_{h}(x^{k}) - \sigma_k s_{g}(y^{k})\\
	& = 2x^{k} - u^{k-1} - \alpha_k s_{h}(x^{k}) - \sigma_k s_{g}(y^{k}) - \tau_{k-1}s_f(u^{k-1}).
	\end{align*}
	Moreover, Step 2 of Algorithm~\ref{a:DY-shadow} and \eqref{eq:OC-v} imply \begin{equation*}
		z^{k+1}  = x^k - \sigma_k s_{g}(y^{k}) - \alpha_k s_h(x^k),
	\end{equation*} which combined with \eqref{eq:OC-u} yields  \begin{align*}
		u^k & = z^{k+1} - \tau_k s_f(u^k)\\
		& = x^k - \alpha_k \left( \frac{\tau_k}{\alpha_k}s_f(u^k) + \frac{\sigma_k}{\alpha_k}s_{g}(y^{k}) + s_h(x^k)\right).
	\end{align*} Setting \begin{equation} \label{eq:dk-shadow}
	d^k = \frac{\tau_k}{\alpha_k}s_f(u^k) + \frac{\sigma_k}{\alpha_k}s_{g}(y^{k}) + s_h(x^k),
	\end{equation} yields $u^k = x^k - \alpha_k d^k$. Thus, from Step 4 of Algorithm~\ref{a:DY-shadow}, it follows that \begin{equation} \label{eq:diff-x} 
	x^{k+1} = x^k - \alpha_k\left(d^k -\frac{\gamma_k}{\alpha_k}(x^k - x^{k-1})\right). 
	\end{equation} 
	Furthermore, subtracting the expressions in Step 4 of Algorithm~\ref{a:DY-shadow} 
at iterations $k$ and $k-1$ yields 
\[u^k - u^{k-1} + \gamma_k(x^k - x^{k-1}) - \gamma_{k-1}(x^{k-1} - x^{k-2}) - x^{k+1} + x^k = 0.\] 
Rearranging terms and using \eqref{eq:diff-x}, we obtain that
 \begin{align*}
		u^k &=u^{k-1} - \gamma_k(x^k - x^{k-1}) + \gamma_{k-1}(x^{k-1} - x^{k-2}) -  \alpha_k\left(d^k - \frac{\gamma_k}{\alpha_k}(x^k - x^{k-1})\right)\\
		&=u^{k-1} -  \alpha_k\left(d^k - \frac{\gamma_{k-1}}{\alpha_k}(x^{k-1} - x^{k-2})   \right).
	\end{align*}

	Since $\{x^k\}$, $\{u^k\}$ and $\{y^k\}$ are bounded by assumption, 
in view of \eqref{eq:dk-shadow} and of the
 local boundedness of the Clarke subdifferential, $\{d^k\}$ is  bounded. In this manner, there exists $c>0$ such that \begin{equation*}
		\max\left\{\|y^k\|, \|d^k\|, \|\diff^k\|\right\} \leq c, \text{~where~} \diff^k = x^k - x^{k-1}.
	\end{equation*} 
Hence, enlarging $X$ if necessary, we guarantee that
 \begin{equation*}
	d^k - \frac{\gamma_k}{\alpha_k}(x^k - x^{k-1}) \in X, \text{~and~} d^k - \frac{\gamma_{k-1}}{\alpha_k}(x^{k-1} - x^{k-2})   \in X.
	\end{equation*} 
Altogether, \begin{equation*}
	\phi_x = d^k - \frac{\gamma_k}{\alpha_k}(x^k - x^{k-1})\quad
  \text{and} \quad \phi_u = d^k - \frac{\gamma_{k-1}}{\alpha_k}(x^{k-1} - x^{k-2}) 
	\end{equation*} satisfy \eqref{eq:def-Phi} with $$d=d^k, \quad \diff^+ = \diff^k, \quad \diff = \diff^{k-1}, \quad x = x^k, \quad u = u^{k-1}, \text{~and~} y = y^k.$$ 
	In particular, observe that the discussion above ensures that for $x,u \in X$, $\PhiD(k,x,u) \subseteq X \times X$, and thus  $\PhiD$ is well-defined.\end{proof}

	\subsection{Convergence of DYS Algorithm~\ref{a:DY-shadow}} \label{ss:convergence-DYS}
	
	In this section, we use the following additional notation. 
Given a set $B \subseteq \R^d$, we denote by $\Delta_B = \{(b,b): b \in B\} $ the diagonal set
in $\R^d\times \R^d$ associated with $B$.
	
	Since Algorithm~\ref{a:DY-shadow} can be seen as an abstract subgradient method 
in view of Proposition~\ref{p:DY-shadow-as-subgradient-step}, its convergence 
properties can be analyzed using Theorem~\ref{t:subgrad-scheme-convergence}. 
 First, in the next lemma, we verify two of the assumptions 
which are required by this theorem.
	
	\begin{lemma} \label{l:aux-results-convergence}
		
		Under the assumptions of Proposition~\ref{p:DY-shadow-as-subgradient-step},  consider the operator $\PhiD$ defined in \eqref{eq:def-Phi}. Then, the following hold: \begin{equation} \label{eq:conv-assumpt-1}
			\sup_{x\in X} \limsup_{\substack{k \to \infty \\ \hat{x} \to x\\ \hat{u} \to x}} \sup_{(\phi_x,\phi_u)\in\PhiD(k,\hat{x},\hat{u})}\|((\phi_x,\phi_u))\| < \infty,
		\end{equation} namely, condition \eqref{t:subgrad-scheme-convergence-1} holds for $\PhiD$. Moreover, for all $x \in X$, 
		 \begin{equation} \label{eq:conv-assumpt-2}
		\tPhiD(x,x) \subseteq \Delta_{\partial F(x) }. 
		\end{equation} 

	\end{lemma}
	
	\begin{proof}
	By the construction above, 
$\tPhiD(k, \cdot,\cdot)$ takes values in $X \times X$ and maps to subsets of $X \times X$. 
Hence, the outer limits of $\Phi(k,x,u)$ are bounded, so that \eqref{eq:conv-assumpt-1} follows.

		To show \eqref{eq:conv-assumpt-2}, 
let $(\phi_x,\phi_u) \in \R^d \times \R^d$ be such that there exist iteration
indices $\ell_k \to \infty$ as $k \to \infty$, and associated (sub)sequences 
 $\hat{x}^k \to x$, $\hat{u}^k \to x$, $(\hat{\phi}_x^k,\hat{\phi}_u^k) \in \PhiD(\ell_k, \hat{x}^k, \hat{u}^k)$ and $(\hat{\phi}_x^k,\hat{\phi}_u^k) \to (\phi_x,\phi_u)$. In this manner, there exist $\widehat{\diff}^k_+, \widehat{\diff}^k \in \R^d$, $\|\widehat{\diff}^k_+\|, \|\widehat{\diff}^k\| \leq c$, $s_h(\hat{x}^k) \in \partial h(\hat{x}^k)$, \[s_g(\hat{y}^k) \in \partial g(\hat{y}^k) \text{~with~} \hat{y}^k = 2\hat{x}^k - \hat{u}^k - \alpha_{\ell_k}s_h(\hat{x}^k) - \sigma_{\ell_k} s_g(\hat{y}^k) - \tau_{\ell_k-1}s_f(\hat{u}^k), \, \|\hat{y}^k\|\leq c, \]  and \[s_f(\hat{u}^k_+) \in \partial f(\hat{u}^k_+), \text{~with~} \hat{u}_+^{k} = \hat{u}^k - \alpha_{\ell_k}\left(\hat{d}^k - \frac{\gamma_{\ell_k -1}}{\alpha_{\ell_k}}\widehat{\diff}^k \right),\] such that for \begin{equation} \label{eq:hat-dk}
			\hat{d}^k = s_h(\hat{x}^k) + \frac{\sigma_{\ell_k}}{\alpha_{\ell_k}}s_g(\hat{y}^k) + \frac{\tau_{\ell_k}}{\alpha_{\ell_k}}s_f(\hat{u}^k_+),
		\end{equation}  it holds that \begin{equation} \label{eq:hat-phik}
		\hat{\phi}^k_x = \hat{d}^k - \frac{\gamma_{\ell_k}}{\alpha_{\ell_k}}\widehat{\diff}^k_+ , \quad \hat{\phi}^k_u =\hat{d}^k - \frac{\gamma_{\ell_k-1}}{\alpha_{\ell_k}}\widehat{\diff}^k.
		\end{equation}  In view of boundedness of the sequences $\{\hat{x}^k\}$, $\{\hat{u}^k\}$ and $\{\hat{y}^k\}$, and local boundedness of the Clarke subdifferential, it follows that $\{\hat{d}^k\}$ is bounded. Furthermore, \eqref{eq:stepsize-assumption} implies $\hat{u}^k_+ \to x$  and $\hat{y}^k \to x$.  Hence, from \eqref{eq:hat-dk}, \eqref{eq:stepsize-assumption} and upper semicontinuity of the Clarke subdifferential, any cluster point of $\{\hat{d}^k\}$ belongs to $\partial h (x) + \partial g(x) + \partial f(x)$. Resorting to \eqref{eq:stepsize-assumption} once again, the fact that both 
$\{\widehat{\diff}^k_+\}$ and $\{\widehat{\diff}^k\}$ are bounded, 
\eqref{eq:hat-phik} and uniqueness of the limit yield that
\[\hat{\phi}^k_x \to \phi_x \in \partial f(x) + \partial g(x) + \partial h(x),\quad
\hat{\phi}^k_u \to \phi_u \in \partial f(x) + \partial g(x) + \partial h(x)\quad\text{and}\quad\phi_x  = \phi_u,\] 
from where we conclude \eqref{eq:conv-assumpt-2}.\end{proof}

The next result states that the shadow sequence $\{u^k\}$ and its inertial version $\{x^k\}$ 
``agree'' asymptotically. This is a direct consequence of $\gamma_k \to 0$ 
and Step 4 in Algorithm~\ref{a:DY-shadow}.

    \begin{lemma}\label{l:same-clusters}
        Consider the sequences $\{x^k\}$ and $\{u^k\}$ generated by Algorithm~\ref{a:DY-shadow}. If assumption \eqref{eq:stepsize-assumption} holds and $\{x^k\}$ is bounded, then $x^k - u^{k-1} \to 0$ and $\acc(\{x^k\}) = \acc(\{u^{k-1}\})$.
    \end{lemma}

In the context of the framework of Section~\ref{sec:abstract} and
of Theorem~\ref{t:subgrad-scheme-convergence}, 
we take the following as the pseudo-Lyapunov function:
	\begin{equation} \label{eq:def-L-DYS}
		\mathcal{L}(x,u) :=\frac{1}{2}\big(F(x) + F(u)\big) .
	\end{equation} By Assumption~\ref{as:DY}\eqref{as:DY-i}, $\mathcal{L}$ is locally Lipschitz continuous and Clarke regular on $\R^d$, so that the hypotheses in Theorem~\ref{t:subgrad-scheme-convergence} on $\mathcal{L}$ hold.

	In the following result, we show that the set of attractors for sequences generated
according to $\PhiD$ defined by \eqref{eq:A}, for a set $C$ to be specified, 
is a subset of the stationary points of problem  \eqref{eq:DY-problem}. In this manner, we relate the DYS abstract subgradient dynamics of Proposition~\ref{p:DY-shadow-as-subgradient-step} to the optimization problem intented to be solved by Algorithm~\ref{a:DY-shadow}.

	\begin{proposition}
		Let Assumption~\ref{as:DY} hold, let $\mathcal{L}$ be defined by
 \eqref{eq:def-L-DYS}, and let $X \subseteq \R^d$ be the set defined in Proposition~\ref{p:DY-shadow-as-subgradient-step}. Let $\{x^k\}$ and $\{u^k\}$ be generated by Algorithm~\ref{a:DY-shadow}. Define the set $Y$ as  \begin{equation} \label{eq:Y-set-DYS}
			Y := \clconv\big\{\{x^k: k \in \N\}\cup\{u^k : k \in \N\}\big\} \subseteq X.
		\end{equation} Let $C = \Delta_Y$, and let $S$ denote the set of stationary points of problem \eqref{eq:DY-problem}. Then \begin{equation} \label{eq:claim-A}
			A(\mathcal{L}, C, \PhiD) \subseteq \Delta_S.
		\end{equation} 
	\end{proposition}

	\begin{proof}
	We shall establish the claim by proving that the negation of \eqref{eq:claim-A}
is false. To that end, let $a \notin \Delta_S$. If $a \notin C$, then 
$a \notin A(\mathcal{L}, C, \PhiD)$ by definition. We can thus assume 
that $a \in C$, which implies the existence of $x \in Y \setminus S$ such that $a = (x,x)$. In this manner, $\bar{s} = \proj_{\partial F(x)}(0) \neq 0$ satisfies \begin{equation} \label{eq:proj-aux}
			\|\bar{s}\|^2 \leq \left<\bar{s}, z \right> \text{~for all~} z \in \partial F(x)
		\end{equation}  and, from \eqref{eq:R}, $\frac{1}{2}(\bar{s},\bar{s}) \in \conv\{\partial \mathcal{L}(x,x)\} \subseteq R(x,x)$. 
Hence, in view of this last inclusion, \eqref{eq:conv-assumpt-2} in 
Lemma~\ref{l:aux-results-convergence} and \eqref{eq:proj-aux}, we obtain that
 \begin{equation*}
			\max_{r \in R(x,x)}\min_{\phi\in\thickbar{\Phi}^{\rm DY}(x,x)} \left<r,\phi\right> \geq \frac{1}{2}\min_{\phi\in\thickbar{\Phi}^{\rm DY}(x,x)} \left<(\bar{s},\bar{s}),\phi\right> \geq \frac{1}{2}\min_{\phi\in \partial F(x)} \left<\bar{s},\phi\right> + \left<\bar{s},\phi\right> \geq \|\bar{s}\|^2>0.
		\end{equation*} 
		It follows from \eqref{eq:A} that $a = (x,x) \notin A(\mathcal{L}, C, \PhiD)$, which proves the claim in \eqref{eq:claim-A}.\end{proof}

We are now in position to establish the main result of this section,
concerning convergence of DYS algorithm.
	
	\begin{theorem}[Convergence of Algorithm~\ref{a:DY-shadow}] \label{t:convergence-DY}

	Let Assumption~\ref{as:DY} hold.	Let $\{x^k\}$, $\{u^k\}$ and $\{y^k\}$ be any bounded sequences generated by Algorithm~\ref{a:DY-shadow}.  Then, the following hold.
        \begin{enumerate}[(i)]
            \item \label{t:convergence-DY-i} There exists $j \in J$ from \eqref{eq:connected-components} such that  \begin{equation*}
            	\acc(\{F(x^k)\}) = \acc(\{F(u^{k-1})\}) \subseteq F\big(S^j \cap \acc(\{x^k\})\big),
            \end{equation*} where $S^j$ is the connected component of $S$ containing $A^j$. 
            \item \label{t:convergence-DY-ii} If, in addition, the set $F(S)$ is nowhere dense in $\R$, then all the accumulation points of $\{x^k\}$ and $\{u^{k-1}\}$ belong to an 
$F$-connected component of $S$.
        \end{enumerate}
	\end{theorem}

	\begin{proof}
		Let $X \subseteq \R^d$ be the open bounded convex
 set from Proposition~\ref{p:DY-shadow-as-subgradient-step}. Let $Y \subseteq \R^d$ be defined in 
\eqref{eq:Y-set-DYS}, a compact convex set containing $\{x^k\}$ and $\{u^k\}$.

		 Furthermore, in view of the
 continuity of $F$ and Lemma~\ref{l:same-clusters}, for any $k_j\to \infty$ such that 
$x^{k_j} \to x^*$ and $u^{k_j-1} \to x^*$, we have that
 $F(x^{k_j}) \to F^*$ if and only if $F(u^{k_j-1}) \to F^*$, which shows that $\acc(\{F(x^k)\}) = \acc(\{F(u^{k-1})\})$. Denote by $K\subseteq \R^d$ the (common) set of accumulation points of $\{x^k\}$ and $\{u^{k-1}\}$.
Then $\Delta_K$ is the set of accumulation points of the sequence $\{(x^k, u^{k-1})\}$.

		To apply Theorem~\ref{t:subgrad-scheme-convergence}, let $V = X\times X$, $v^k = (x^k, u^{k-1})$, $\Phi$ defined in \eqref{eq:def-Phi}, and recall that $C = \Delta_Y \subseteq X \times X$. Observe that \eqref{eq:gen-bounded-assumption} holds directly, condition \eqref{t:subgrad-scheme-convergence-1} holds from \eqref{eq:conv-assumpt-1} in Lemma~\ref{l:aux-results-convergence}, and \eqref{t:subgrad-scheme-convergence-2} is  part of the stepsizes assumptions. Then, Theorem~\ref{t:subgrad-scheme-convergence} applies and the following hold. 
		
		First, there exists $j \in J$ such that  all accumulation points of the sequence $\{\mathcal{L}(x^k,u^{k-1})\}$ belong to the set $\mathcal{L}(A^j \cap \Delta_K)$.   
 From \eqref{eq:claim-A}, $A^j \subseteq \Delta_S$, and since $F\big(\pi_1(A^j)\big)=\mathcal{L}(A^j)$ is connected, it follows that $\pi_1(A^j)$ is an $F$-connected subset $S$, and thus belongs to an $F$-connected component of $S$, denoted $S^j$. 
Observe that \begin{equation} \label{eq:claim-1}
		    \mathcal{L}(A^j \cap \Delta_K) \subseteq F(S^j \cap K).
		\end{equation}  
		Indeed, let $\ell \in \mathcal{L}(A^j \cap \Delta_K)$. Then there exists $(x,u) \in A^j \cap \Delta_K$ such that $\ell = \mathcal{L}(x,u)$. In particular, $x = u \in S$ in view of \eqref{eq:claim-A}, $x\in\pi_1(A^j) \cap K$, and $\ell = F(x) $. 
		This proves our claim. 

Next, let $F^*$ be an accumulation point of $\{F(x^k)\}$. After taking a subsequence if necessary, there exists $k_j \to \infty$ such that $F(x^{k_j}) \to F^*$, $F(u^{k_j-1}) \to F^*$,  $u^{k_j-1} \to x^*$ and $x^{k_j} \to x^*$. Then, from the
 continuity of $F$ and uniqueness of the limit, $\mathcal{L}(x^{k_j},u^{k_j-1}) = \frac{1}{2}(F(x^{k_j}) + F(u^{k_j-1})) \to F(x^*) = F^*$. Therefore, from the first conclusion from Theorem~\ref{t:subgrad-scheme-convergence} and \eqref{eq:claim-1}, $F^* \in \mathcal{L}(A^j \cap \Delta_K) \subseteq F(S^j \cap K)$, which proves \eqref{t:convergence-DY-i}.
		
		Suppose now that $F(S)$ is nowhere dense in $\R$. The relation in \eqref{eq:claim-A} implies $\mathcal{L}(A(\mathcal{L}, C, \Phi)) \subseteq \mathcal{L}(\Delta_S) = F(S)$, hence 
		the second consequence of Theorem~\ref{t:subgrad-scheme-convergence} is that the set $\Delta_K$ of accumulation points of $\{(x^k, u^{k-1})\}$ is a subset of an $\mathcal{L}$-connected component of $A(\mathcal{L}, C, \Phi)$. In view of \eqref{eq:claim-A}, $\Delta_K$ is also a subset of an $\mathcal{L}$-connected component of $\Delta_S$, denoted $W$.  Then $K = \pi_1(\Delta_K) \subseteq \pi_1(W) \subseteq S$. 
        Furthermore, $(F \circ \pi_1)(W) = \mathcal{L}(W)$ is connected in $\R$, where the equality follows from the fact that $W \subseteq \Delta_S$. Since $\{A^j\}_{j \in J}$ is a partition of $A(\mathcal{L}, C, \Phi)$ into $\mathcal{L}$-connected components, then for any $j \in J$, $A^j$ is a maximal subset in $\R^d \times \R^d$ with the property that $\mathcal{L}(A^j)$ is connected in $\R$. As $A^j\subseteq \Delta_S$ in view of \eqref{eq:claim-A}, maximality is also with respect to subsets of $\Delta_S$. Moreover, observe that $A^j\subseteq \Delta_S$ if and only  $A^j = \Delta_{U^j}$ for $U^j := \pi_1(A^j)\subseteq S$. Then, the following holds true since $\pi_1$ is an injection on $\Delta_S$: \begin{align*}
				&A^j \text{~is maximal among subsets of $\Delta_S$ such that~} \mathcal{L}(A^j) \text{~is connected}\\
				\iff& \text{there is no $\tilde{A} \subseteq \Delta_S$, $A^j\subsetneq \tilde{A}$, such that   } \mathcal{L}(\tilde{A})  = F(\pi_1(\tilde{A})) \text{~is connected}\\
				\iff & \text{there is no $\tilde{U} \subseteq S$, $\pi_1(A^j)\subsetneq \tilde{U}$, such that   } F(\tilde{U}) \text{~is connected}\\
				\iff & \pi_1(A^j)  \text{~is maximal among subsets of $S$ such that~} F(\pi_1(A^j)) \text{~is connected}.
		\end{align*} In this manner, $K$ is a subset of the $F$-connected set $\pi_1(W) \subseteq S$ and, in particular, a subset of an $F$-connected component of $S$. This shows \eqref{t:convergence-DY-ii} and concludes the proof.\end{proof}

    \begin{remark}
        In addition to being finite locally Lipschitz and Clarke regular functions, if $f$, $g$ and $h$ are definable in a o-minimal structure (e.g., semialgebraic), then $F(S)$ is finite \cite[Corollary 9~(ii)]{bolte2007clarke}, and thus, in particular, nowhere dense in $\R$. Hence, the additional condition in Theorem~\ref{t:convergence-DY} \eqref{t:convergence-DY-ii} holds.
    \end{remark}

	In this section, we proved that for the DYS method to work in nonconvex optimization, it is enough to assume local Lipschitz continuity and Clarke regularity of the model components, and no smoothness condition is required to compensate nonconvexity. Furthermore,  our approach allows us to include momentum to a shadow sequence of DYS in a natural way.

	\section{Nonconvex PDHG method without smoothness or strong convexity} \label{s:PDHG}

    In this section, we investigate convergence properties of a generalization of PDHG method
\cite{chambolle2011first,condat2013primal,vu2013splitting}
tailored to solve the following optimization problem \begin{equation} \label{eq:CP-problem}
		\min_{x \in \R^d} f(Lx) + g(x) + h(x),
	\end{equation} where $f: \R^p \to \R$, $g, \, h: \R^d \to \R$ can
 be nonsmooth nonconvex, and $L: \R^d \to \R^p$ is a linear operator. The version of PDHG we analyze is given in Algorithm~\ref{a:PDHG}, 
which results from replacing in \eqref{eq:PDHG} the gradient step on $h$ with a subgradient step, 
and the proximal step on $f$ (after using the Moreau decomposition that relates the proximal 
operators of $f$ and $f^*$ in the convex setting) and the proximal step on $g$, 
with the corresponding stationary proximal steps. 
Our goal is to show that this method finds a stationary point of problem \eqref{eq:CP-problem}, namely, a point $\bar{x} \in \R^d$ such that \begin{equation*}
	0 \in L^* \partial f(L\bar{x})  + \partial g(\bar{x}) + \partial h(\bar{x}).\end{equation*} 

    \begin{algorithm}[!ht]
	\caption{PDHG for solving problem \eqref{eq:CP-problem}. \label{a:PDHG}}
	\SetKwInOut{Input}{Input}
	\Input{Choose starting points $x^0 \in \R^d$, $z^1 \in \R^p$.  
		Set $x^1 = u^0 = x^0$. }
	
	\For{$k=0,1,2,\dots$}{
		Step 1. Choose $\sigma_k > 0$ and compute $y^{k+1} \in \argstat_{\sigma_k^{-1}f}(Lu^k + \sigma_k^{-1}z^k)$. 
		
		Step 2. Update $z^{k+1} = z^k + \sigma_k(Lu^k - y^{k+1})$.
		
		Step 3.  Choose $\tau_k, \, \alpha_k >0$, and compute $s_{h}(x^{k}) \in \partial h(x^k) $ and $$x^{k+1} \in \argstat_{\tau_k g}(x^k - \alpha_k s_h(x^k) - \alpha_k L^* z^{k+1}).$$
		
		Step 4. Choose $\theta_{k+1}>0$ and update $u^{k+1} = x^{k+1} + \theta_{k+1}(x^{k+1} - x^{k})$.
	}
\end{algorithm}

    Next, we specify the nature of the functions $f,\, g$ and $h$, and the parameters of Algorithm~\ref{a:PDHG}.

    \begin{assumption} \label{as:PDHG}
		For Algorithm~\ref{a:PDHG}, suppose the following conditions hold.
		\begin{enumerate}[(i)]
			\item \label{as:PDHG-C-i} The functions $f: \R^p \to \R$, $g: \R^d \to \R$ and $h : \R^d \to \R$ are locally Lipschitz continuous and Clarke regular, and $L : \R^d \to \R^p$ is a (bounded) linear map.
			\item \label{as:PDHG-C-ii} The parameters of Algorithm~\ref{a:PDHG} satisfy
			\begin{equation} \label{as:stepsize-PDHG}
				\frac{\tau_k}{\alpha_k} \to 1, \quad \alpha_k \to 0, \quad \sum_{k=1}^\infty \alpha_k = \infty, \quad \sigma_k \to \infty, \quad \text{and} \quad \bar{\theta} := \sup_{k} \theta_k < \infty.
			\end{equation}
			
		\end{enumerate}
	\end{assumption}

    \begin{remark}\label{r:regime-PDHG}
    Two comments on the parameters assumption \eqref{as:stepsize-PDHG} are in order.
    \begin{enumerate}
            \item[(i)] The assumption $\sigma_k \to \infty$ makes Step~1 of Algorithm~\ref{a:PDHG} well-posed for nonconvex $f$. For example, if $f$ is $r_f$-weakly convex, the mapping $\argstat_{\sigma_k^{-1}f}$ is single-valued and Lipschitz continuous as soon as $\sigma_k > r_f$ \cite[Proposition 12.19]{rockafellar2009variational}.
            
            \item[(ii)] Since $\tau_k \sim \alpha_k \to 0$, the product $\tau_k\sigma_k$ can still be kept below the classical PDHG threshold. Indeed, suppose $h$ is differentiable with $\ell_h$-Lipschitz continuous gradient. Then, taking $\tau_k = \alpha_k$ and $\sigma_k = \rho/\alpha_k$ with $\rho>0$ such that $\rho \|L\|^2 \leq 1$, yields \begin{equation*}
            	\dfrac{\tau_k \ell_h}{2} + \tau_k\sigma_k\|L\|^2 = \dfrac{\tau_k \ell_h}{2} + \frac{\tau_k}{\alpha_k}\rho\|L\|^2 \to \rho\|L\|^2  \leq 1.
            \end{equation*}
            
        \end{enumerate}
	\end{remark}

    In this section we show that the primal sequence $\{x^k\}$ generated by Algorithm~\ref{a:PDHG}
 can be viewed as a generalized subgradient sequence in the sense of Section~\ref{sec:abstract},
 for the objective
	\begin{equation} \label{eq:F-PDHG}
		F := f \circ L + g + h,
	\end{equation} and that its cluster points, if any, are stationary for \eqref{eq:CP-problem}. Throughout this section we set
	\begin{equation*}
		S := \{ x \in \R^d : 0 \in \partial F(x)\}.
	\end{equation*} Observe that due to Clarke regularity, $\partial F(x) = L^*\partial f\circ (Lx) + \partial g(x) + \partial h(x) $ for all $x \in \R^d$.
	
	\subsection{The PDHG abstract subgradient scheme} \label{ss:scheme}
	
	Let us record elementary consequences of Assumption~\ref{as:PDHG} that will be used repeatedly. 

    \begin{lemma}
        Consider the sequences $\{y^k\}$, $\{z^k\}$ and $\{u^k\}$ generated by Algorithm~\ref{a:PDHG}. Then, the following hold.
        \begin{enumerate}
            \item For all $k \in \N$,
	\begin{equation} \label{eq:z-is-sf}
		z^{k+1} = s_f(y^{k+1}) \in \partial f(y^{k+1}), \qquad \|y^{k+1} - Lu^k\| = \frac{\|z^k - z^{k+1}\|}{\sigma_k}.
	\end{equation} 

            \item If, in addition, the sequences $\{x^k\}$ and $\{z^k\}$ generated by Algorithm~\ref{a:PDHG} are bounded,  the constants \begin{equation}\label{eq:constants-PDHG}
            \begin{aligned}
		&M_z := \sup_k \|z^k\|, \quad M_x := \sup_k \|x^k\|, \quad M_y := \sup_k \|y^k\|, \\ &M_\phi := \sup_k \left\| s_h(x^k) + \frac{\tau_k}{\alpha_k}s_g(x^{k+1}) + L^*z^{k+1} \right\|
        \end{aligned}
	\end{equation} satisfy \begin{equation} \label{eq:c-PDHG}
		c := \max\{M_x, M_y, M_z, M_\phi\} < \infty.
	\end{equation} 
    \end{enumerate}
    \end{lemma}
    \begin{proof}
    From Steps~1 and~2 of Algorithm~\ref{a:PDHG}, there exists  $s_f(y^{k+1}) \in \partial f(y^{k+1}) $ such that \begin{equation*}
        0  = s_f(y^{k+1}) - z^k + \sigma_k(y^{k+1} - Lu^k) =  s_f(y^{k+1}) - z^{k+1},
    \end{equation*} and thus \eqref{eq:z-is-sf} follows. Moreover, $M_z < \infty$ and $M_x < \infty$ are direct from the assumption. Furthermore, since $\bar{\theta}<\infty$ from \eqref{as:stepsize-PDHG}, the sequence $\{u^k\}$ is bounded. By \eqref{eq:z-is-sf} and $\sigma_k \to \infty$, so is $\{y^k\}$. Finally, in view of local boundedness of the Clarke subdifferential and $\tau_k/\alpha_k \to 1$, $M_\phi$ is finite.\end{proof}

	In the following, we interpret Algorithm~\ref{a:PDHG} as an instance of the abstract 
subgradient scheme of Section~\ref{sec:abstract}.
	
\begin{proposition} \label{p:scheme-PDHG} Let Assumption~\ref{as:PDHG} hold. Suppose the sequences $\{x^k\}$ and $\{z^k\}$ generated by Algorithm~\ref{a:PDHG} are bounded. Consider an open 
	convex bounded set $X \subseteq \R^d$ containing $\{x^k\}$.  For a constant $c >0$ dependent on $X$, let the set-valued mapping $\PhiP: \N \times X  \to 2^X $ be given by, for all $x \in \R^d $, \begin{equation} \label{eq:def-Phi-PDHG}
		\phi \in \PhiP(k,x) \iff\left\{\begin{aligned}
			& \phi = s_h(x) + \frac{\tau_k}{\alpha_k} s_g(x_+) + L^* s_f(y_+), \\
			& u = x + \theta_k \delta, \quad \|\delta\| \leq c\,\alpha_{k-1},  \\
			& s_h(x) \in \partial h(x),\, s_g(x_+) \in \partial g(x_+),\, s_f(y_+) \in \partial f(y_+), \\
			& 0 = s_f(y_+) - z + \sigma_k(y_+ - Lu), \\
			& 0 = \tau_k s_g(x_+) + x_+ - (x - \alpha_k s_h(x)) + \alpha_k L^* s_f(y_+),  \\
			& \|y_+\| \leq c, \quad \|z\| \leq c, \quad \|x_+\| \leq c.
		\end{aligned}\right.
	\end{equation}
	Then,  \begin{equation} \label{eq:scheme-PDHG}
		x^{k+1} \in x^k - \alpha_k \PhiP(k,x^k) \quad \text{for all } k \in \N.
	\end{equation} 
\end{proposition}
\begin{proof}
	Fix $k \in \N$. From Steps~1 and~3 in Algorithm~\ref{a:PDHG}, there exist $s_f(y^{k+1}) \in \partial f(y^{k+1})$ and $s_g(x^{k+1}) \in \partial g(x^{k+1})$ such that \begin{align}
		0 & = s_f(y^{k+1}) - z^k + \sigma_k(y^{k+1} - Lu^k), \label{eq:step-1-PDHG} \\
		0 & = \tau_k s_g(x^{k+1}) + x^{k+1} - x^k + \alpha_k s_h(x^k) + \alpha_k L^* z^{k+1}. \label{eq:step-3-PDHG} 
	\end{align} Set
	\begin{equation} \label{eq:phi-k-PDHG}
		\phi^k := s_h(x^k) + \frac{\tau_k}{\alpha_k}s_g(x^{k+1}) + L^* z^{k+1},
	\end{equation} and recall from \eqref{eq:z-is-sf} that $s_f(y^{k+1}) = z^{k+1}$. Choose, in \eqref{eq:def-Phi-PDHG},
	\begin{equation*}
		\delta = x^k - x^{k-1}, \quad z = z^k, \quad y_+ = y^{k+1}, \text{~and}\quad x_+ = x^{k+1}.
	\end{equation*} 
	Note that the  fourth line of \eqref{eq:def-Phi-PDHG} is \eqref{eq:step-1-PDHG}, and the   fifth line is \eqref{eq:step-3-PDHG} after substituting $z^{k+1} = s_f(y^{k+1})$. Rearranging \eqref{eq:step-3-PDHG} and using \eqref{eq:phi-k-PDHG} yields
	\begin{equation*}
		x^{k+1} = x^k - \alpha_k\left( s_h(x^k) + \frac{\tau_k}{\alpha_k}s_g(x^{k+1}) + L^*s_f(y^{k+1}) \right) = x^k - \alpha_k \phi^k.
	\end{equation*} From Step~4 of Algorithm~\ref{a:PDHG} we have $u = u^k = x^k + \theta_k(x^k-x^{k-1}) = x + \theta_k\delta $. For $k = 1$, we have $\delta = 0$ by the convention $x^0 = x^1$. For $k \geq 2$, the update $x^k = x^{k-1} - \alpha_{k-1}\phi^{k-1}$, 
	together with \eqref{eq:constants-PDHG} and \eqref{eq:c-PDHG}, give
	\begin{equation*}
		\|\delta\| = \alpha_{k-1}\|\phi^{k-1}\| \leq \alpha_{k-1}M_\phi \leq \alpha_{k-1}\,c.
	\end{equation*} The bounds for $\|z\| = \|z^k\|$, $\|y_+\| = \|y^{k+1}\|$ and $\|x_+\| = \|x^{k+1}\|$ are immediate from \eqref{eq:constants-PDHG} and \eqref{eq:c-PDHG}. Hence $\phi^k \in \PhiP(k,x^k)$ and \eqref{eq:scheme-PDHG} holds. Observe that, enlarging $X$ if necessary, local boundesness of the Clarke subdifferential ensures that $\PhiP(k,x) \subseteq X$ whenever $x \in X$.\end{proof}

\begin{remark}
	Observe that the dual variable and the momentum term are absorbed into the generalized subgradient set-valued map as auxiliary quantities.
\end{remark}

\subsection{Convergence of PDHG Algorithm~\ref{a:PDHG}} \label{ss:convergence-PDHG}

We now proceed to show that the interpretaton of Algorithm~\ref{a:PDHG} 
in Proposition~\ref{p:scheme-PDHG} complies with the general abstract subgradient framework
of Section~\ref{sec:abstract}. 

\begin{lemma} \label{l:lim-Phi-PDHG}
	Let Assumption~\ref{as:PDHG} hold, and consider the operator $\PhiP$ defined in \eqref{eq:def-Phi-PDHG}. Then,
	\begin{equation*} 
		\sup_{\hat{x} \in X} \limsup_{\substack{x \to \hat{x} \\ k \to \infty}} \sup_{\phi \in \PhiP(k,x)} \|\phi\| < \infty,
	\end{equation*} that is, condition \eqref{t:subgrad-scheme-convergence-1} holds for $\PhiP$. Furthermore, for all $x \in X$, \begin{equation*}
		\tPhiP(x) \subseteq \partial h(x) + \partial g(x) + L^* \partial f(Lx).
	\end{equation*}
\end{lemma}

\begin{proof}
	
	Let $x \in X$, $k \in \N$ and $\phi \in \PhiP(k,x)$. Using the notation in 
	Proposition~\ref{p:scheme-PDHG}, since $X$ is bounded and $\|x_+\| \leq c$, $\|y_+\|\leq c$, 
	local boundedness of the Clarke subdifferential implies that
	\begin{equation*}
		\|\phi\| \leq \sup_{w \in \overline{X}} \sup_{s \in \partial h(w)}\|s\| + \left(\sup_k \frac{\tau_k}{\alpha_k}\right)\sup_{\|w\|\leq c} \sup_{s \in \partial g(w)}\|s\| + \|L\| \sup_{\|w\|\leq c} \sup_{s \in \partial f(w)}\|s\|,
	\end{equation*} where the right-hand side is finite because $\tau_k/\alpha_k \to 1$. As this bound is independent of $x$ and $k$, the first claim follows.

	Let $x \in X$ and $\phi \in \tPhiP(x)$. Then there exist $\ell_k \to \infty$, $\hat{x}^k \to x$ and $\hat{\phi}^k \to \phi$ with $\hat{\phi}^k \in \Phi(\ell_k, \hat{x}^k)$ for all $k \in \N$. Let $\hat{\delta}^k$, $\hat{z}^k$, $\hat{y}^k_+$, $\hat{x}^k_+$ and the corresponding subgradients conform to
	\eqref{eq:def-Phi-PDHG} at $(\ell_k, \hat{x}^k)$, and set $\hat{u}^k = \hat{x}^k + \theta_{\ell_k}\hat{\delta}^k$. Denote
	\begin{equation*}
		M := \sup_{\|w\| \leq c} \sup_{s \in \partial f(w)}\|s\| < \infty.
	\end{equation*}  By \eqref{eq:def-Phi-PDHG}, $\|\hat{\delta}^k\| \leq c\,\alpha_{\ell_k - 1}$. Since $\alpha_j \to 0$ and $\ell_k \to \infty$, we have $\alpha_{\ell_k-1} \to 0$, and since $\theta_{\ell_k} \leq \bar{\theta}$,
	\begin{equation*}
		\|\hat{u}^k - \hat{x}^k\| \leq \bar{\theta}c\,\alpha_{\ell_k-1} \to 0 .
	\end{equation*} Together with $\hat{x}^k \to x$, this gives $\hat{u}^k \to x$ and, 
	by the  continuity of $L$,
	$L\hat{u}^k \to Lx$. Next, from the fourth line of \eqref{eq:def-Phi-PDHG}, $\|\hat{z}^k\| \leq c$, $\|\hat{y}^k_+\|\leq c$ and $\sigma_j \to \infty$, it follows that
	\begin{equation*}
		\|\hat{y}^k_+ - L \hat{u}^k\| = \frac{\|\hat{z}^k - s_f(\hat{y}^k_+)\|}{\sigma_{\ell_k}} \leq \frac{c+M}{\sigma_{\ell_k}} \to 0,
	\end{equation*} and thus $\hat{y}^k_+ \to Lx$. Furthermore, since $\{\hat{x}^k_+\}$, $\{\hat{x}^k\}$ and $\{\hat{y}^k_+\}$ are bounded, local boundedness of the Clarke subdifferential implies that $\{s_g(\hat{x}^k_+)\}$, $\{s_h(\hat{x}^k)\}$ and $\{s_f(\hat{y}^k_+)\}$ are bounded as well. The fifth line of \eqref{eq:def-Phi-PDHG} can be written as
	\begin{equation*}
		\hat{x}^k_+ = \hat{x}^k - \tau_{\ell_k}s_g(\hat{x}^k_+) - \alpha_{\ell_k}\left( s_h(\hat{x}^k) + L^* s_f(\hat{y}^k_+)\right).
	\end{equation*}  As $\tau_j \to 0$ and $\alpha_j \to 0$ by \eqref{as:stepsize-PDHG}, we conclude that
	$\hat{x}^k_+ - \hat{x}^k \to 0$, hence $\hat{x}^k_+ \to x$. Up to a subsequence, we have $s_h(\hat{x}^k) \to \zeta_h$, $s_g(\hat{x}^k_+) \to \zeta_g$ and $s_f(\hat{y}^k_+) \to \zeta_f$. Outer semicontinuity of the Clarke subdifferential gives $\zeta_h \in \partial h(x)$, $\zeta_g \in \partial g(x)$ and $\zeta_f \in \partial f(Lx)$. From using $\tau_{\ell_k}/\alpha_{\ell_k} \to 1$ in the first line of \eqref{eq:def-Phi-PDHG}, it follows that
	\begin{equation*}
		\hat{\phi}^k = s_h(\hat{x}^k) + \frac{\tau_{\ell_k}}{\alpha_{\ell_k}}s_g(\hat{x}^k_+) + L^* s_f(\hat{y}^k_+) \longrightarrow \zeta_h + \zeta_g + L^*\zeta_f.
	\end{equation*} As $\hat{\phi}^k \to \phi$, uniqueness of the limit yields $\phi = \zeta_h + \zeta_g + L^*\zeta_f \in \partial h(x) + \partial g(x) + L^* \partial f(Lx)$, which proves the claim.\end{proof}

To apply Theorem~\ref{t:subgrad-scheme-convergence}, here we take as the
pseudo-Lyapunov function the objective function itself, namely,
\begin{equation} \label{eq:def-L-PDHG}
	\mathcal{L} := F = f\circ L + g + h .
\end{equation} By Assumption~\ref{as:PDHG}\eqref{as:PDHG-C-i}, $\mathcal{L}$ is locally Lipschitz continuous and Clarke regular on $\R^d$, so that the standing hypotheses of Theorem~\ref{t:subgrad-scheme-convergence} on $\mathcal{L}$ are satisfied.

We next prove that the set of attractors of the iterates
generated according to $\PhiP$ given by \eqref{eq:A} is a subset 
of the set of stationary points of problem \eqref{eq:CP-problem}.

\begin{proposition} \label{p:claim-A-PDHG}
	Let Assumption~\ref{as:PDHG} hold, let $\mathcal{L}$ be given by \eqref{eq:def-L-PDHG}, and let $X \subseteq \R^d$ be the set defined in Proposition~\ref{p:scheme-PDHG}. Let $\{x^k\}$ be the sequence generated by Algorithm~\ref{a:PDHG}. Define the set $C$ as \begin{equation} \label{eq:C-set-PDHG}
		C := \clconv\{x^k : k \in \N\} \subseteq X.
	\end{equation}  Then
	\begin{equation} \label{eq:claim-A-PDHG}
		A(\mathcal{L},C,\PhiP) \subseteq  S.
	\end{equation}
\end{proposition}

\begin{proof}
	Let $x \notin S$, that is, $0 \notin \partial F(x)$. If $x \notin C$, then $x \notin A(\mathcal{L},C,\PhiP)$ from definition. Suppose then that $x \in C$. We show that $x \notin A(\mathcal{L},C,\PhiP)$. Since the set $\partial F(x)$ is nonempty, convex and compact, it follows that  the projection $\bar{r} := \proj_{\partial F(x)}(0)$ is well-defined, and $\bar{r} \neq 0$. The characterization of the projection onto a closed convex set gives $\left<0-\bar{r}, w - \bar{r}\right> \leq 0$ for all $w \in \partial F(x)$, that is, $\left<\bar{r},w\right> \geq \|\bar{r}\|^2$ for all $w \in \partial F(x)$. Note that $\bar{r} \in \partial F(x) = \partial \mathcal{L}(x) \subseteq R(x)$, and $\tPhiP(x) \subseteq \partial F(x)$ by Lemma~\ref{l:lim-Phi-PDHG}. Then
	\begin{equation*}
		\max_{r \in R(x)} \min_{\phi \in \tPhiP(x)} \left<r,\phi\right> \geq \min_{\phi \in \tPhiP(x)} \left<\bar{r},\phi\right> \geq \min_{w \in \partial F(x)}\left<\bar{r},w\right> \geq \|\bar{r}\|^2 > 0.
	\end{equation*} In view of \eqref{eq:A}, $x \notin A(\mathcal{L},C,\PhiP)$, concluding the proof.\end{proof}

We conclude this section by establishing our
convergence result for Algorithm~\ref{a:PDHG}.

\begin{theorem}[Convergence of Algorithm~\ref{a:PDHG}] \label{t:convergence-PDHG}
	Let Assumption~\ref{as:PDHG} hold. 
	Let the sequences $\{x^k\}$ and $\{z^k\}$  generated by Algorithm~\ref{a:PDHG} be bounded. 
	Then the following hold.
	\begin{enumerate}[(i)]
		\item \label{t:convergence-C-i} There exists $j \in J$ from \eqref{eq:connected-components} such that all the accumulation points of the sequence $\{F(x^k)\}$ belong to $F\big(A^j \cap \acc(\{x^k\})\big)$. 
		
		\item \label{t:convergence-C-ii} If, in addition, the set $F(S)$ is nowhere dense in $\R$, then $\acc(\{x^k\})$ is contained in a single $F$-connected component of $A(\mathcal{L},C,\PhiP)$. In particular, every accumulation point of $\{x^k\}$ is a stationary point of $F$.  
		
		Furthermore, $\acc(\{x^k\})  = \acc(\{u^k\})  $, $\acc(\{y^k\}) \subseteq L(\acc(\{x^k\}) ) $, and for any $\bar x \in \acc(\{x^k\})$, there exists $ \bar z \in \acc(\{z^k\})$ such that \begin{equation*}
			0 \in L^* \bar z + \partial g(\bar x) + \partial h(\bar x) \quad \text{and} \quad \bar z \in \partial f (L \bar x).
		\end{equation*}
	\end{enumerate}
\end{theorem}

\begin{proof}
	
	Let $X \subseteq \R^d$ be the set defined in Proposition~\ref{p:scheme-PDHG}, and $C$ given in \eqref{eq:C-set-PDHG}. Take $V = X$ and $v^k = x^k$ for all $k \in \N$. Consider $\mathcal{L}$ given by \eqref{eq:def-L-PDHG}, and $\PhiP$ given by \eqref{eq:def-Phi-PDHG}. Observe that condition \eqref{eq:gen-bounded-assumption} holds by assumption and the definition of $C$, which is compact and convex and contains $\acc(\{x^k\})$. The scheme \eqref{eq:scheme-PDHG} and condition \eqref{t:subgrad-scheme-convergence-1} hold by Proposition~\ref{p:scheme-PDHG}, and \eqref{t:subgrad-scheme-convergence-2} is contained in \eqref{as:stepsize-PDHG}. Since $\mathcal{L}(x^k) = F(x^k)$, Theorem~\ref{t:subgrad-scheme-convergence} yields the existence of  $j \in J$ from \eqref{eq:connected-components} such that all the accumulation points of the sequence $\{F(x^k)\}$ belong to $F\big(A^j \cap \acc(\{x^k\})\big)$, which gives item~\eqref{t:convergence-C-i}. 
	
	To prove \eqref{t:convergence-C-ii}, as we have 
	$\mathcal{L}(A(\mathcal{L},C,\PhiP)) = F(A(\mathcal{L},C,\PhiP))$, 
	Proposition~\ref{p:claim-A-PDHG} ensures that
	$ \mathcal{L}(A(\mathcal{L},C,\PhiP)) \subseteq F(S)$. Since a subset of a nowhere dense set is nowhere dense, the set $\mathcal{L}(A(\mathcal{L},C,\PhiP))$ is nowhere dense in $\R$, and the second part of Theorem~\ref{t:subgrad-scheme-convergence} implies that $\acc(\{x^k\})$ is a subset of an $\mathcal{L}$-connected component of $A(\mathcal{L},C,\PhiP)$. Invoking \eqref{eq:claim-A-PDHG} once more, this component is contained in $S$. 
	This shows the first part of \eqref{t:convergence-C-ii}.

	Let $\bar x$ be an accumulation point of $\{x^k\}$, and take a subsequence $\{x^{k_j}\}$ such that $x^{k_j} \to \bar x \in S$ as $j \to \infty$. Observe that $\|x^{k+1}-x^k\| \to 0$ from \eqref{eq:scheme-PDHG}, $\alpha_k \to 0$, and $\PhiP(k,x^k) \subseteq X$.  In view of Step 4 
	in Algorithm~\ref{a:PDHG}, it follows that $u^{k_j} \to \bar x $. 
	Moreover, from \eqref{eq:constants-PDHG}, since $\{z^k\}$ is bounded and $\sigma_k \to \infty$, it follows that $y^{k+1} - Lu^k \to 0$. Hence, $y^{k_j+1} \to L \bar x$. Up to a subsequence, we also have that $z^{k_j + 1} \to \bar z$ for some $\bar z \in \R^p$. Using \eqref{eq:constants-PDHG} again and the fact that $\partial f$ is outer semicontinuous, $\bar z \in \partial f (L \bar x)$. This concludes the proof.\end{proof}

\begin{remark}
	
	Some comments about Theorem~\ref{t:convergence-PDHG} are in order.
	
	\begin{enumerate}
		\item If $f$ is globally $\ell_f$-Lipschitz continuous, then $\|z^{k+1}\| = \|s_f(y^{k+1})\| \leq \ell_f$ for every $k$ by \eqref{eq:z-is-sf}. Hence, we only need boundedness of the primal sequence $\{x^k\}$.
		
		\item In the context of Theorem~\ref{t:convergence-PDHG}\eqref{t:convergence-C-ii}, if $f$ is convex we recover a KKT point of problem	\eqref{eq:CP-problem}. Indeed, the pair $(\bar x, \bar z)$ solves \begin{equation*}
			0 \in L^* \bar z + \partial g(\bar x) + \partial h(\bar x) \quad \text{and} \quad  L \bar x \in \partial f^* (\bar z),
		\end{equation*} the system of equations that defines the stationary points of the following min-max formulation \begin{equation*}
			\min_{x \in \R^d}\max_{z \in \R^p} \; g(x) + h(x) + \langle Lx , z \rangle - f^*(z).
		\end{equation*} This results generalizes, for example,  
		the original convex case with smooth $h$ \cite{condat2013primal}, and the 
		weakly convex -- strongly convex setting with $h=0$ in \cite{mollenhoff2015primal}.
		
	\end{enumerate}
\end{remark}

\section*{Declarations}

\noindent {\small
{\bf Conflicts of Interest.}
The authors declare that they have no conflict of interest of any kind related to the manuscript.}

\noindent {\small
{\bf Data Availability Statement.}
Data sharing is not applicable to this article as no datasets were generated or analyzed
hereby.}

\noindent {\small
{\bf Artificial Intelligence Statement.} 
The authors declare that they have not used AI tools.}

\noindent {\small
{\bf Funding.} The first author was partially supported by Centro de Modelamiento Matem\'atico (CMM), ACE210010
and FB210005, BASAL funds for centers of excellence from ANID-Chile. The second author is supported 
in part by CNPq Grant 306775/2023-9 and
by FAPERJ Grant E-26/200.165/2026.}

\printbibliography
\end{document}